\documentclass[preprint,12pt]{elsarticle}

\usepackage{amssymb}
\usepackage{amsmath}
\usepackage{zymacros1}
\usepackage{makecell}
\usepackage{multirow,makecell}
\usepackage{tablefootnote}
\usepackage{placeins}
\usepackage{epstopdf}
\usepackage{subcaption}

\usepackage{amsthm}
\newtheorem{theorem}{Theorem}%  meant for continuous numbers
\newtheorem{proposition}[theorem]{Proposition}
\newtheorem{lemma}[theorem]{Lemma}
\newtheorem{corollary}[theorem]{Corollary}%
\newtheorem{remark}{Remark}
\newtheorem{assumption}{Assumption}%

\newtheorem{definition}{Definition}

\begin{document}

\begin{frontmatter}



\title{Nearly Optimal Approximation Rates for Deep Super ReLU Networks on Sobolev Spaces} %% Article title

\author[label1]{Yahong Yang$\footnote{Corresponding author(\textit{yyang3194@gatech.edu})}$}
\affiliation[label1]{organization={School of Mathematics, Georgia Institute of Technology},
            addressline={686 Cherry Street},
            city={Atlanta},
            postcode={30332},
            state={Georgia},
            country={USA}}

\affiliation[label4]{organization={Department of Mathematics and Department of Computer Science, University of Maryland},
            addressline={4176 Campus Dr.},
            city={College Park},
            postcode={20742},
            state={Maryland},
            country={USA}}

\affiliation[label2]{organization={Department of Mathematics},
            addressline={The Hong Kong University of Science and Technology, Clear Water Bay},
            city={Kowloon},
            postcode={999077},
            state={Hong Kong SAR},
            country={China}}

\affiliation[label3]{organization={Algorithms of Machine Learning and Autonomous Driving Research Lab},
            addressline={HKUST Shenzhen-Hong Kong Collaborative Innovation Research Institute},
            city={Shenzhen},
            postcode={518048},
            state={Guangdong},
            country={China}}

\author[label2,label3]{Yue Wu}\author[label4]{Haizhao Yang} \author[label2,label3]{Yang Xiang}%% Author name

% \affiliation{organization={},%Department and Organization

%% Abstract
\begin{abstract}
This paper introduces deep super ReLU networks (DSRNs) as a method for approximating functions in Sobolev spaces measured by Sobolev norms $W^{m,p}$ for $m\in\mathbb{N}$ with $m\ge 2$ and $1\le p\le +\infty$. Standard ReLU deep neural networks (ReLU DNNs) cannot achieve this goal. DSRNs consist primarily of ReLU DNNs, and several layers of the square of ReLU added at the end to smooth the networks output. This approach retains the advantages of ReLU DNNs, leading to the straightforward training. The paper also proves the optimality of DSRNs by estimating the Vapnik–Chervonenkis dimension of higher-order derivatives of DNNs, and obtains the generalization error in Sobolev spaces via an estimate of the pseudo dimension of higher-order derivatives of DNNs.
\end{abstract}


% \begin{highlights}
% \item We introduce deep super ReLU networks for approximating functions in Sobolev spaces, measured by norms in $W^{m,p}$ for $m \ge 2$. In our DNNs, ReLU is the predominant activation function, with the square of ReLU appearing only in the final layers. This approach enhances the smoothness of the DNNs while maintaining a simple network structure without intricate composition formulas.
% \item We estimate the upper bound of the VC-dimension of higher-order derivatives of DNNs. By utilizing our estimation of the VC-dimension, we demonstrate the optimality of our DNN approximation, as measured by Sobolev norms.
% \item We provide an upper bound estimation for the pseudo-dimension of high-order derivatives of DNNs. By utilizing this estimation, we can obtain the nearly optimal generalization error of DNNs in Sobolev spaces that contain high-order derivatives of DNNs.
% \end{highlights}

%% Keywords

\end{frontmatter}



\section{Introduction}
 Deep neural networks (DNNs) with the rectified linear unit (ReLU) activation function \cite{glorot2011deep} have become increasingly popular in scientific and engineering applications, including image classification \cite{krizhevsky2017imagenet,he2015delving}, regularization \cite{czarnecki2017sobolev}, and dynamic programming \cite{finlay2018lipschitz,werbos1992approximate}. However, when tasked with approximating higher-order derivatives of target functions, such as in solving partial differential equations (PDEs) \cite{Lagaris1998,weinan2017deep,raissi2019physics,weinan2018deep} by DNNs, ReLU DNNs are unable to provide accurate results due to their lack of smoothness. For instance, physics-informed neural networks (PINNs) \cite{raissi2019physics} have emerged as a popular approach to solving PDEs. When solving a PDE on the domain $\Omega = [a,b]^d$,
\begin{equation}
\begin{cases}
\fL u = f & \text{in } \Omega, \\
u = g & \text{on } \partial\Omega,
\end{cases} \notag
\end{equation}
the loss function $L(\vtheta)$ is defined as
\[
L(\vtheta)
:= \|\fL \phi(\vx; \vtheta) - f\|_{L^2(\Omega)}^2
+ \beta \, \|\phi(\vx; \vtheta) - g\|_{L^2(\partial\Omega)}^2,
\]
where $\phi(\vx; \vtheta)$ denotes a deep neural network (DNN) with parameters $\vtheta$,  
$\beta > 0$ is a weight balancing the boundary and interior terms, and $\fL$ is a second-order differential operator.  
However, ReLU-based DNNs are unsuitable for this task because they lack the smoothness required to compute $\fL \phi(\vx; \vtheta)$ accurately; in fact, $\fL \phi(\vx; \vtheta)$ may not even exist in the classical sense.

In spite of the inherent smoothness limitations of ReLU-based deep neural networks, they are still widely used in training due to their simplicity, efficient computation, scale-invariance, and the fact that they are less prone to vanishing gradients. On the other hand, some smooth activation functions, such as sigmoid functions, can suffer from vanishing gradients, making them unsuitable for specific tasks. Additionally, certain smooth activation functions, such as the square of ReLU, can have complex composition formulas making training difficult.

In this paper, we propose a method for using DNNs to approximate functions in Sobolev space measured by the norms in Sobolev spaces. Specifically, the DNNs we construct in this paper use two types of activation functions, ReLU and the square of ReLU. The square of ReLU is only present in the last several layers of the DNNs, and the number of layers with the square of ReLU is proportional to the order of $\log_2 L$, where $L$ is the depth of the DNNs. We call this kind of DNNs as deep super ReLU networks (DSRNs).

The approximation by neural networks in various spaces measured by variable norms has been established in a lot of papers. The following table (Table~\ref{tab:summary}) lists the references and shows the difference between their work and ours. \begin{table}[h!]
\centering
\setlength{\tabcolsep}{6pt}
\renewcommand{\arraystretch}{1.2}
\begin{tabular}{|c|c|c|c|}
\hline
Measured Norm / Target Space &  Activations & {Architecture } &  References \\
\hline
$L^p$/Sobolev spaces
  & G.A\tablefootnote{G.A (general activations) indicates a broad activation class: either multiple specific activations are treated, or a generalized activation satisfying stated conditions.}
  & \multirow{2}{*}{Shallow}
  & \cite{mhaskar1996neural,liu2025integral,yang2024nonparametric} \\
\cline{1-2}\cline{4-4}
\makecell{$L^2$, $L^\infty$, $H^1$\\/(spectral) Barron spaces}
  & G.A
  & 
  & \cite{weinan2019barron,siegel2022sharp,siegel2023characterization,klusowski2018approximation} \\
\hline
$L^p$/Smooth spaces
  & ReLU
  & \multirow{7}{*}{Deep}
  & \cite{shen2019nonlinear,shen2020deep,lu2021priori,petersen2018optimal} \\
\cline{1-2}\cline{4-4}
$L^p$/Sobolev spaces
  & ReLU and \(\tanh\)
  &
  & \cite{yarotsky2017error,siegel2022optimal,de2021approximation} \\
\cline{1-2}\cline{4-4}
\makecell{$L^p$ and $H^1$\\/Special spaces}
  & \makecell{ReLU or RePU}
  &
  & \makecell{\cite{montanelli2019new,suzukiadaptivity,blanchard2021shallow,yang2024near}\\ \cite{opschoor2024exponential,opschoor2022exponential}}\\
\cline{1-2}\cline{4-4}
$L^p$/Smooth spaces
  & G.A
  &
  & \cite{zhang2024deep} \\
\cline{1-2}\cline{4-4}
\makecell{$W^{\alpha,p}$, \(0\le\alpha\le 1\)\\/Sobolev spaces}
  & ReLU
  &
  & \cite{yang2023nearly,guhring2020error} \\
\cline{1-2}\cline{4-4}
$W^{m,\infty}$/$W^{n,\infty}$
  & G.A
  & rate $\epsilon^{-\frac{d}{n-m}}$
  & \cite{guhring2021approximation,he2023deep,he2023expressivity} \\
\cline{1-2}\cline{4-4}
$W^{m,p}$/$W^{n,p}$
  & \makecell{ReLU and\\squared ReLU}
  & rate $\epsilon^{-\frac{d}{2(n-m)}}$
  & \textit{\makecell{Theorem~\ref{main2} and\\ Corollary~\ref{main4}}} \\
\hline
\end{tabular}
\caption{Summary of approximation results across various norms/spaces, activation choices, and architectures for fully connected neural networks. In the last two rows, the rate indicates the number of parameters required for a neural network to achieve an $\epsilon$ approximation accuracy measured in the specified norms.}
\label{tab:summary}
\end{table}
The above table summarizes only a portion of the approximation results for fully connected neural networks. Other network architectures, such as convolutional neural networks~\cite{zhou2020universality,feng2021generalization,mao2022approximation} and transformers~\cite{cole2025context,havrilla2024understanding}, are not discussed here, as they are not closely related to the focus of our work.

In \cite{yang2023nearly}, ReLU and its squared version were utilized as DNN activation functions, yet the following three questions remain. 

Firstly, the problem of approximating Sobolev functions with DNNs exhibiting super convergence rates, measured by $W^{m,p}$ norms for $p\in[1,\infty)$, remains unresolved, although the $W^{m,\infty}$ case was addressed by \cite{yang2023nearly}, and the case of $m=0$ was considered by \cite{siegel2022optimal}. Secondly, proper distributions of the squared ReLU function remain unknown. Incorporating it at any layer of a DNN with a depth of $L$ might lead to a polynomial degree of $2^L$ with respect to parameters, which poses challenges in training. Thirdly, the optimality of the approximation process remains unclear due to the complexity associated with higher-order derivatives of DNNs using ReLU and squared ReLU activations. To establish optimality, estimating the optimal Vapnik-Chervonenkis dimension (VC-dimension) and pseudo-dimension \cite{anthony1999neural, vlassis2021sobolev, abu1989vapnik, pollard1990empirical} of their higher order derivatives, which characterize the complexity or richness of a function set, need to be estimated. Even for DNNs without derivatives, as demonstrated in \cite{bartlett2019nearly}, obtaining optimal bounds is challenging. The upper VC-dimension and pseudo-dimension for DNNs using squared ReLU is $\fO(N^2L^3\log_2 N\log_2L)$ or $\fO(N^3L^2\log_2 N\log_2L)$. However, this result is not directly applicable to proving the optimality of squared ReLU-based DNN approximation due to a polynomial gap in depth $L$.

\begin{figure}[h!]
\centering
\includegraphics[scale=0.57]{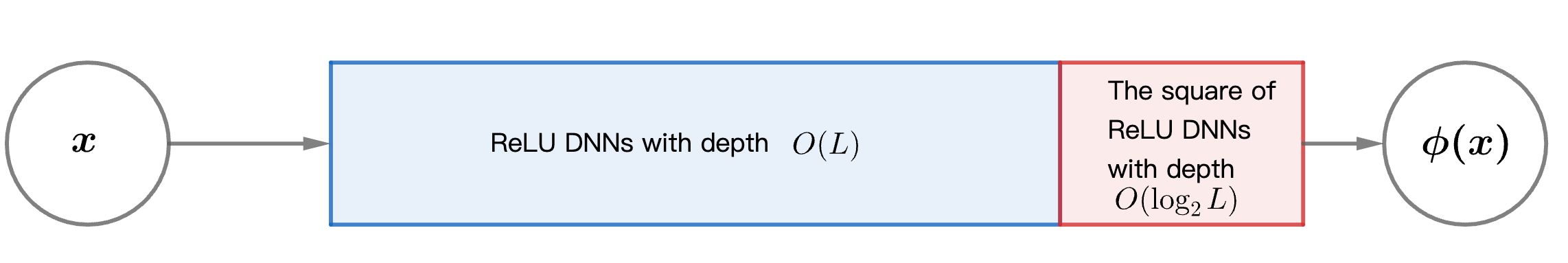}
\caption{The structure of DSRNs.}
\label{DSRN}
\end{figure}

In this paper, the deep super ReLU networks (DSRNs) is presented to address these three questions, as illustrated in Fig.~\ref{DSRN}. The DSRN architecture consists of two components for a network with depth $L$: ReLU DNNs with a depth of approximately $\fO(L)$, and the square of ReLU DNNs with a depth of approximately $\fO(\log_2L)$—significantly shallower than the first component.

The key achievement of this paper is that we demonstrate this type of DNNs can effectively approximate functions within Sobolev spaces as measured by Sobolev norms $W^{m,p}$ for $m\ge 2$. This result offers benefits both in theory and training practice. From a theoretical perspective, the proposed approach yields a nearly optimal approximation rate concerning the width $N$ and depth $L$ of DSRNs. This is made possible by the fact that the square of ReLU only appears in a shallower layer, approximately $\fO(\log_2L)$, thereby maintaining VC-dimension and pseudo-dimension within $\fO(N^2L^2\log_2 N\log_2L)$ and eliminating the polynomial gap in the squared ReLU-based DNNs. {In \cite{yang2023nearly}, the optimality of deep neural networks was proven only for networks where all activation functions are ReLU. The inclusion of mixed activation functions increased the complexity of the neural network space, creating a gap between space complexity and the approximation rate due to the frequent occurrence of the square of ReLU.}

In terms of training, our DSRNs predominantly employ ReLU activations, with the square of ReLU appearing only at the end of the DNNs. This enhances the smoothness of DNNs while preserving a straightforward network structure without intricate composition formulas. {This contrasts with \cite{yang2023nearly}, which did not specify where to apply the square of ReLU, leading to training difficulties. Our experiments in Sec.~\ref{num} demonstrate that our proposed DNN (DSRN) trains effectively, unlike those using the square of ReLU throughout, which are difficult to train due to gradient vanishing. For the training of DSRN, we discuss more in Section~\ref{training} later.

Finally, we estimate the generalization error of DSRNs in Theorem~\ref{general thm} under Sobolev training \cite{czarnecki2017sobolev}.  
Although there are works focusing on the generalization error for Sobolev training—such as for physics-informed neural networks \cite{lorenz2024error,jiang2024generalization,mishra2022estimates,de2024error} and the Deep Ritz method \cite{jiao2021error,lu2021priori,yang2023nearly,jiao2025drm}—most existing results consider continuous approximators.  
In these cases, the approximation rate is consistent with classical results \cite{devore1989optimal,devore1993constructive}, where the parameter size is uniformly bounded.  
Under this setting, the generalization error can be bounded by directly estimating the covering number of the parameter space.

However, in the generalization analysis of deep neural networks, if we aim to achieve super-convergence—as shown in Theorem~\ref{main3}—many parameters can become large and are not uniformly bounded.  
Therefore, a different approach is required to estimate the generalization error in Sobolev training when the network achieves super-convergence.  
Our results are nearly optimal with respect to the number of sample points and significantly improve upon related works such as \cite{yang2023nearly,jiao2021error}.  
In Theorem~\ref{general thm}, the obtained generalization rate is nearly optimal with respect to the number of sampling points. There are two key reasons for this improvement. First, we employ covering number arguments instead of Rademacher complexity to bound the generalization error. Rademacher complexity is a global measure that may lose local structural information, which can result in suboptimal error rates. To achieve optimal generalization error bounds, it is necessary to adopt localized complexity measures, such as local Rademacher complexity \cite{bartlett2005local,yang2024nonparametric} or covering number techniques \cite{gyorfi2002distribution,liu2024learning}. Both approaches can yield optimal bounds, in this work, we adopt the covering number approach. The second reason lies in our tight estimates of the VC-dimension and pseudo-dimension of higher-order derivatives of deep neural networks. These estimates lead to sharper covering number bounds. To the best of our knowledge, this is the first work that explicitly analyzes the VC- and pseudo-dimensions of higher-order derivatives in neural networks, and applies them directly to the error analysis in Sobolev training.
}

\subsection{Training of DSRN}\label{training}

In this section, we discuss the training of DSRNs, which naturally leads to the following fundamental questions. A DSRN, defined as the composition of a squared ReLU neural network with another ReLU neural network, is generally non-smooth for most parameter choices. This raises two important questions:

\textbf{Questions:} 
(i) Can such a neural network approximate a target function measured in a Sobolev norm?  
(ii) If so, can this type of neural network be effectively trained in practice?

\paragraph{Answer to Question (i) and proof sketch.}
Theorem~\ref{main2} together with Corollary~\ref{main4} gives an affirmative answer. We briefly outline the construction.

\emph{(1) Local polynomial approximation.}
By the Bramble–Hilbert lemma \cite{brenner2008mathematical}, on each cell of a shape-regular partition of the domain there exists a polynomial (of the prescribed order) that approximates the target function in the relevant Sobolev norm with the standard rate.

\emph{(2) Exact realization of polynomials by squared-ReLU.}
Squared-ReLU networks can represent polynomials exactly; more specifically, a degree-$n$ polynomial can be realized by a squared-ReLU network of depth $\mathcal{O}(\log_2 n)$ \cite{he2023expressivity}. Hence, on each cell we implement the local polynomial by a very shallow squared-ReLU subnetwork (serving as a smooth “feature extractor”).

\emph{(3) Gating and global assembly.}
To “turn on’’ the appropriate local feature, we use a (standard) ReLU network to implement piecewise-constant gating weights; an efficient construction follows the bit-extraction paradigm \cite{bartlett2019nearly}. This yields a global piecewise model. The only nonsmoothness arises along the interfaces where neighboring local pieces interact (the “bad regions’’).

\emph{(4) Smoothing via partitions of unity within the architecture.}
We build several such subnetworks using differently shifted local partitions so that their interface sets do not coincide. We then combine them with a smooth partition of unity (PoU) whose supports are chosen to avoid the bad region of at least one subnetwork on each point of the domain. Because DSRNs can implement $(x)_+^{s}$ for any $s\in\mathbb{N}$, they can represent truncated-power bases and hence B-splines; consequently, smooth PoUs (and smooth piecewise-polynomial functions) can be realized within the DSRN \cite{de1978practical}. This yields the desired Sobolev-norm approximation claimed in Theorem~\ref{main2} and Corollary~\ref{main4}.

\begin{figure}[h!]
    \centering
    \includegraphics[scale=0.09]{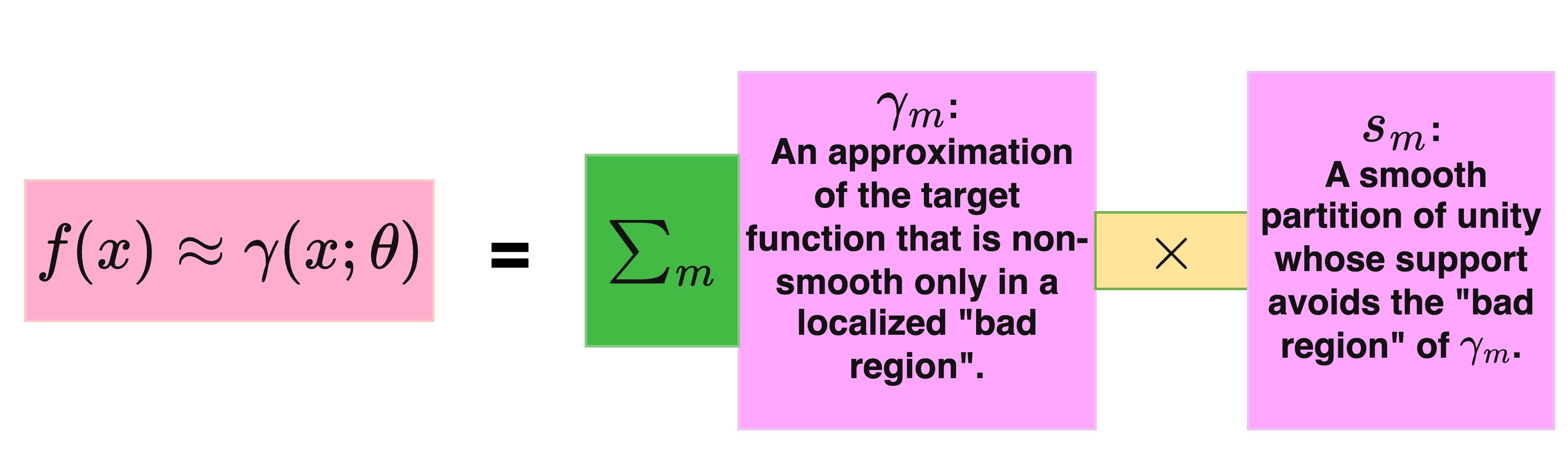}
    \caption{The architecture of the smooth DSRN.}
    \label{instr}
\end{figure}

Further mathematical details can be found in Section~\ref{sketchproof}. The overall structure of the neural network is illustrated in Fig.~\ref{instr}.

For the second question, we divide the discussion into two parts.  

First, can we actually run the code? The answer is yes. Although the initial guess of a neural network may be non-smooth, the non-smooth points lie on a set of measure zero in the domain, which is unlikely to be sampled. Even if some points are sampled at these locations, we can simply resample during training to avoid them. Therefore, the experiments are indeed feasible.  

Second, can such neural networks learn functions and match higher-order derivatives? The answer is also affirmative. During training, the DSRN gradually becomes smooth. After random initialization, the DSRN outputs exhibits jump discontinuities, leading to highly irregular derivatives. However, since the loss function involves high-order derivatives, the training process gradually reduces these discontinuities. 
The theoretical analysis guarantees that there exist globally smooth DSRNs capable of approximating target functions in Sobolev norms. Consequently, the network systematically eliminates jump discontinuities and converges toward smooth function approximations with respect to higher-order Sobolev norms.

In practice, there may still remain very small jump points, but their effect is negligible, and the network can already approximate the derivatives of the target function well. If even higher smoothness is desired, one can resample around these jump points to further reduce the gaps.  

\begin{figure}[!ht]
\centering
\includegraphics[scale=0.35]{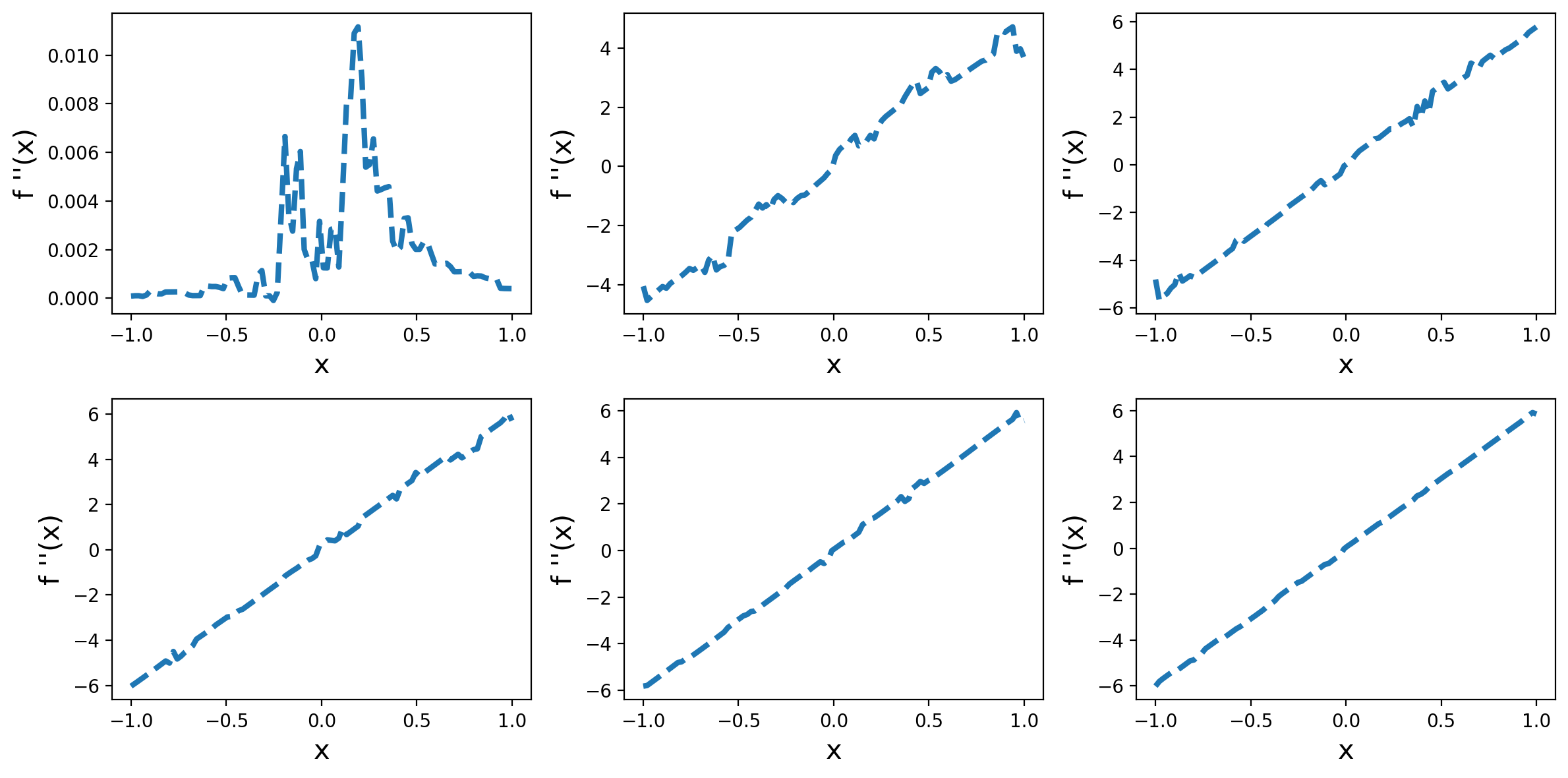}
\caption{Evolution of the second-order derivative estimated by DSRN during training. The panels (arranged from top to bottom and left to right) correspond to the initial state, iterations 10,000, 20,000, 30,000, 40,000, and the final iteration.}
\label{fig:dynamic 2nd derivative}
\end{figure}

To illustrate this dynamic smoothing process, we use a DSRN to approximate $f(x)=x^3$ 
and
Fig.~\ref{fig:dynamic 2nd derivative} shows the evolution of the second-order derivative estimated by DSRN during the training. The learned second-order derivatives quickly approach the ground truth $f''(x)=6x$ across the domain, demonstrating that DSRN progressively captures the smoothness of the target function in terms of the Sobolev norm. More detailed experiments for solving PDEs can be found in Sec.~\ref{num}.

\subsection{Main contributions}
Our main contributions are:

$\bullet$ We introduce deep super ReLU networks for approximating functions in Sobolev spaces, measured by norms in $W^{m,p}$ for $m \ge 2$. In our DNNs, ReLU is the predominant activation function, with the square of ReLU appearing only in the final layers. This approach enhances the smoothness of the DNNs while maintaining a simple network structure without intricate composition formulas.

$\bullet$ We estimate the upper bound of the VC-dimension of higher-order derivatives of DNNs. By utilizing our estimation of the VC-dimension, we demonstrate the optimality of our DNN approximation, as measured by Sobolev norms.

$\bullet$ We provide an upper bound estimation for the pseudo-dimension of high-order derivatives of DNNs. By utilizing this estimation, we can obtain the nearly optimal generalization error of DNNs in Sobolev spaces that contain high-order derivatives of DNNs.

\section{Preliminaries}\label{preliminaries}
		\subsection{Neural networks}
		Let us summarize all basic notations used in the DNNs as follows:
		
		\textbf{1}. Matrices are denoted by bold uppercase letters. For example, $\vA\in\sR^{m\times n}$ is a real matrix of size $m\times n$ and $\vA^\T$ denotes the transpose of $\vA$.
		
		\textbf{2}. Vectors are denoted by bold lowercase letters. For example, $\vv\in\sR^n$ is a column vector of size $n$. Furthermore, denote $\vv(i)$ as the $i$-th elements of $\vv$.
		
		\textbf{3}. For a $d$-dimensional multi-index $\valpha=[\alpha_1,\alpha_2,\cdots\alpha_d]\in\sN^d$, we denote several related notations as follows: 
  \begin{align}
  &(a)~ |\boldsymbol{\alpha}|=\left|\alpha_1\right|+\left|\alpha_2\right|+\cdots+\left|\alpha_d\right|; \notag\\&(b)~\boldsymbol{x}^{\boldsymbol{\alpha}}=x_1^{\alpha_1} x_2^{\alpha_2} \cdots x_d^{\alpha_d},~ \boldsymbol{x}=\left[x_1, x_2, \cdots, x_d\right]^\T;\notag\\& (c)~\boldsymbol{\alpha} !=\alpha_{1} ! \alpha_{2} ! \cdots \alpha_{d} !.\notag
  \end{align}
		
		%\textbf{4}. Let $B_{r,|\cdot|}(\vx)\subset\sR^d$ be the closed ball with a center $\vx\in\sR^d$ and a radius $r$ measured by the Euclidean distance. Similarly, $B_{r,\|\cdot\|_{\ell_\infty}}(\vx)\subset\sR^d$ be the closed ball with a center $\vx\in\sR^d$ and a radius $r$ measured by the $\ell_\infty$-norm.
		
\textbf{4}. Assume $\vn\in\sN_+^m$, and $f$ and $g$ are functions defined on $\sN_+^m$, then $f(\vn)=\fO(g(\vn))$ means that there exists positive $C$ independent of $\vn,f,g$ such that $f(\vn)\le Cg(\vn)$ when all entries of $\vn$ go to $+\infty$.
		
\textbf{5}. Define $\sigma_1(x):=\max\{0,x\}$ and $\sigma_2:=\sigma_1^2(x)$. We call the neural networks with activation function $\sigma_t$ with $t= i$ as $\sigma_i$ neural networks ($\sigma_i$-NNs), $i=1,2$. With the abuse of notations, we define $\sigma_i:\sR^d\to\sR^d$ as $\sigma_i(\vx)=\left[\begin{array}{c}
			\sigma_i(x_1) \\
			\vdots \\
			 \sigma_i(x_d)
		\end{array}\right]$ for any $\vx=\left[x_1, \cdots, x_d\right]^T \in\sR^d$.
		
\textbf{6}. Define $L,N\in\sN_+$, $N_0=d$ and $N_{L+1}=1$, $N_i\in\sN_+$ for $i=1,2,\ldots,L$, then a $\sigma_i$-NN $\phi$ with the width $N$ and depth $L$ can be described as follows:\[\boldsymbol{x}=\tilde{\boldsymbol{h}}_0 \stackrel{W_1, b_1}{\longrightarrow} \boldsymbol{h}_1 \stackrel{\sigma_i}{\longrightarrow} \tilde{\boldsymbol{h}}_1 \ldots \stackrel{W_L, b_L}{\longrightarrow} \boldsymbol{h}_L \stackrel{\sigma_i}{\longrightarrow} \tilde{\boldsymbol{h}}_L \stackrel{W_{L+1}, b_{L+1}}{\longrightarrow} \phi(\boldsymbol{x})=\boldsymbol{h}_{L+1},\] where $\vW_i\in\sR^{N_i\times N_{i-1}}$ and $\vb_i\in\sR^{N_i}$ are the weight matrix and the bias vector in the $i$-th linear transform in $\phi$, respectively, i.e., $\boldsymbol{h}_i:=\boldsymbol{W}_i \tilde{\boldsymbol{h}}_{i-1}+\boldsymbol{b}_i, ~\text { for } i=1, \ldots, L+1$ and $\tilde{\boldsymbol{h}}_i=\sigma_i\left(\boldsymbol{h}_i\right),\text{ for }i=1, \ldots, L.$ In this paper, an DNN with the width $N$ and depth $L$, means
		(a) The maximum width of this DNN for all hidden layers is less than or equal to $N$.
		(b) The number of hidden layers of this DNN is less than or equal to $L$.

\subsection{Sobolev spaces}
		
Denote $\Omega$ as $[a,b]^d$, $D$ as the weak derivative of a single variable function, $D^{\valpha}=D^{\alpha_1}_1D^{\alpha_2}_2\ldots D^{\alpha_d}_d$ as the partial derivative of a multivariable function, where $\valpha=[\alpha_{1},\alpha_{2},\ldots,\alpha_d]^T$ and $D_i$ is the derivative in the $i$-th variable.
		
		\begin{definition}[Sobolev Spaces \cite{evans2022partial}]
			Let $n\in\sN$ and $1\le p\le \infty$. Then we define Sobolev spaces\[W^{n, p}(\Omega):=\left\{f \in L^p(\Omega): D^{\valpha} f \in L^p(\Omega) \text { for all } \boldsymbol{\alpha} \in \sN^d \text { with }|\boldsymbol{\alpha}| \leq n\right\}\] with a norm $\|f\|_{W^{n, p}(\Omega)}:=\left(\sum_{0 \leq|\alpha| \leq n}\left\|D^{\valpha} f\right\|_{L^p(\Omega)}^p\right)^{1 / p}$, if $p<\infty$, and \[\|f\|_{W^{n, \infty}(\Omega)}:=\max_{0 \leq|\alpha| \leq n}\left\|D^{\valpha} f\right\|_{L^\infty(\Omega)}.\]
			Furthermore, for $\vf=(f_1,f_2,\ldots,f_d)$, $\vf\in W^{1,\infty}(\Omega,\sR^d)$ if and only if $ f_i\in W^{1,\infty}(\Omega)$ for each $i=1,2,\ldots,d$ and \[\|\vf\|_{W^{1,\infty}(\Omega,\sR^d)}:=\max_{i=1,\ldots,d}\{\|f_i\|_{W^{1,\infty}(\Omega)}\}.\]
		\end{definition}
		
		\begin{definition}[Sobolev semi-norm \cite{evans2022partial}]
			Let $n\in\sN_+$ and $1\le p\le \infty$. Then we define Sobolev semi-norm $|f|_{W^{n, p}(\Omega)}:=\left(\sum_{|\alpha|= n}\left\|D^{\valpha} f\right\|_{L^p(\Omega)}^p\right)^{1 / p}$, if $p<\infty$, and $|f|_{W^{n, \infty}(\Omega)}:=\max_{|\alpha| = n}\left\|D^{\valpha} f\right\|_{L^\infty(\Omega)}$. Furthermore, for $\vf\in W^{1,\infty}(\Omega,\sR^d)$, we define  \[|\vf|_{W^{1,\infty}(\Omega,\sR^d)}:=\max_{i=1,\ldots,d}\{|f_i|_{W^{1,\infty}(\Omega)}\}.\]
		\end{definition}

\section{Deep Super ReLU Networks}
This paper focuses on deep neural networks, which are characterized by a large value of $L$. It is well-known that ReLU neural networks are unable to approximate functions in the $W^{m,\infty}$ norm for $m\ge 2$. An innovative neural network architecture is introduced in \cite{yang2023nearly}, which cleverly incorporates both ReLU and squared ReLU activation functions. This architecture proves to be highly effective in the approximation of functions within Sobolev spaces, as evidenced by their alignment with Sobolev norms $W^{m,\infty}((0,1)^d)$ for $m\ge 2$: \begin{corollary}[\cite{yang2023nearly}]\label{main3}
Let \(f\in W^{n,\infty}((0,1)^{d})\) with \(\|f\|_{W^{n,\infty}((0,1)^{d})}\le 1\).
For any \(N,L,m\in\mathbb{N}_{+}\) satisfying
\(NL+2^{\lfloor\log_{2}N\rfloor}\ge\max\{d,n\}\) and
\(L\ge\lceil\log_{2}N\rceil\),
there exists a neural network \(\varphi(\vx)\) that mixes the ReLU and squared-ReLU activations, whose width is \(\mathcal{O}(N\log N)\) and whose depth is \(\mathcal{O}(L\log L)\), such that
\[
\|f-\varphi\|_{W^{m,\infty}((0,1)^{d})}
=\mathcal{O}\!\bigl(N^{-2(n-m)/d}L^{-2(n-m)/d}\bigr).
\]
		\end{corollary}

{Nevertheless, a limitation of this corollary that uses ReLU and its square as activation functions in DNNs is that the distribution of the square of ReLU is yet to be determined. If the distribution is entirely unknown, training complexity may significantly increase. Hence, this paper introduces a novel DNN architecture that predominantly uses ReLU activation functions, allowing the square of ReLU to appear only in the last several layers. This design allows us to achieve a nearly optimal rate of approximation, reduce the complexity of DNNs as stated in Corollary \ref{main3}, and obtain a better generalization error, making the training process easier. We will rigorously discuss this structure later and show by the numerical experiment that this network performs better in training compared to networks where the activation functions are either all ReLU or all square of ReLU in Sec.~\ref{num}.
}

We define a subset of $\sigma_2$-NNs with $L\gg 1$ and $C=\fO(1)$ with respect to $L$ as follows:
\begin{align}
&\fN_{C,L}:=\{\psi(\vx)=\psi_2\circ\vpsi_1 (\vx): \psi_2\text{ is a $\sigma_2$-NN with depth $L_2$,}\notag\\&\text{ each component of $\vpsi_1$ is a $\sigma_1$-NN with depth $L_1$ },L_1+L_2\le L,~ L_2\le C \log L.\}\notag
\end{align} 
We refer to elements in $\fN_{C,L}$ as deep super ReLU networks (DSRNs).

\section{Approximation by DSRNs for Functions in Sobolev Spaces}
In this section, we investigate the ability of deep super ReLU networks (DSRNs) to approximate functions in Sobolev spaces. Specifically, we consider the Sobolev space $W^{m,p}([a,b]^d)$ for $m\ge 2$. We establish that DSRNs can effectively approximate functions in $W^{m,p}([a,b]^d)$ with a nearly optimal rate of convergence. {Without loss of generality, we consider the case of the $m=2$.}

\begin{theorem}\label{main2}
	For any $f\in W^{n,p}([a,b]^d)$ for $p\in[1,\infty]$, any $n\ge2$ and $d,N, L\in\sN_+$ with $\log N\ge d(\log n+\log d)$, there is a DSRN $\gamma(\vx)$ in $\fN_{1,9(L+1)\log_24L}$ with the width $2^{d+6}n^{d+1}(N+d)\log_2(8N)$ such that\[\|f(\vx)-\gamma(\vx)\|_{W^{2,p}([a,b]^d)}\le 2^{d+7}C_3(n,d,a,b)\|f(\vx)\|_{W^{n,p}([a,b]^d)}N^{-2(n-2)/d}L^{-2(n-2)/d},\]where $C_3=C_1+4^{\frac{2d+p}{2p}}C_2$ is the constant independent of $N,L$, and $C_1,C_2$ is defined in Proposition \ref{fN}.
\end{theorem}

{\begin{remark}
    As shown in Theorem~\ref{main2}, the approximation rate of the DSRN is 
    \[
        \fO(N^{-2(n-2)/d}L^{-2(n-2)/d})
    \] 
    using \(\fO(N^2L)\) parameters, which outperforms traditional methods and achieves \textit{super-convergence}. This rate is derived via the bit extraction technique, where the approximation error decays faster than the metric entropy of the Sobolev function class. Consequently, observing this rate in numerical experiments can be difficult, as discussed in detail in \cite{siegel2022optimal}. Furthermore, the constant \(C_3\) appearing in Theorem~\ref{main2}, which determines the leading coefficient of the approximation error, depends on the interval endpoints \(a\) and \(b\). Specifically, \(C_1\) depends on \(\max\{|b - a|^{n - i}\}_{i = 0}^2\), and \(C_2\) depends on \(\max\{|a|, |b|\}\). Further details are provided in Remark~\ref{C}.

\end{remark}

\subsection{Proof Sketch for Theorem \ref{main2}}\label{sketchproof}
The proof of Theorem \ref{main2} can be outlined in five parts, and the complete proof is provided in next subsections:
		
		\textbf{(i)}: First, define a sequence of subsets of $\Omega$ in Definition \ref{omega}.

        Then we define a partition of unity $\{g_{\vm}\}_{\vm\in\{1,2\}^d}$ on $[a,b]^d$ with $\operatorname{supp}g_{\vm}\cap[a,b]^d\subset \Omega_{\vm}$ for each $\vm\in\{1,2\}^d$ in Definition \ref{sm}.
		
		\textbf{(ii)}: Next we use the Proposition \ref{smm} to present $\{s_{\vm}\}_{\vm\in\{1,2\}^d}$ by $\sigma_2$-NNs.
		
		\textbf{(iii)}: Third, for each $\Omega_{\vm}\subset [a,b]^d$, where $\vm\in\{1,2\}^d$, we find a function $f_{K,\vm}$ satisfying \begin{align}
			\|f-f_{K,\vm}\|_{W^{2,p}(\Omega_{\vm})}&\le C_1(n,d,|b-a|)\|f\|_{W^{n,p}(\Omega_{\vm})}K^{-(n-2)},\notag\\\|f-f_{K,\vm}\|_{W^{1,p}(\Omega_{\vm})}&\le C_1(n,d,|b-a|)\|f\|_{W^{n,p}(\Omega_{\vm})}K^{-(n-1)},\notag\\\|f-f_{K,\vm}\|_{L^{p}(\Omega_{\vm})}&\le C_1(n,d,|b-a|)\|f\|_{W^{n,p}(\Omega_{\vm})}K^{-n}, 
		\end{align} where $C_1$ is a constant independent of $K$. Moreover, each $f_{K,\vm}$ can be expressed as $f_{K,\vm}=\sum_{|\valpha|\le n-1}g_{f,\valpha,\vm}(\vx)\vx^{\valpha}$, where $g_{f,\valpha,\vm}(\vx)$ is a piecewise constant function on $\Omega_{\vm}$. The proof of this result is based on the Bramble-Hilbert Lemma \cite[Lemma 4.3.8]{brenner2008mathematical}.
		
		\textbf{(iv)}: The fourth step involves approximating $f_{K,\vm}$ using DSRN $\gamma_{\vm}$, following the approach outlined in \cite{lu2021deep}. This method is suitable for our work because $g_{f,\valpha,\vm}(\vx)$ is a piecewise constant function on $\Omega_{\vm}$, and the weak derivative of $g_{f,\valpha,\vm}(\vx)$ on $\Omega_{\vm}$ is zero. This property allows for the use of the $L^\infty$ norm approximation method presented in \cite{lu2021deep}. Thus, we obtain a neural network $\psi_{\vm}$ with width $\fO(N\log N)$ and depth $\fO(L\log L)$ such that \begin{align}\|f_{K,\vm}-\gamma_{\vm}(\vx)\|_{W^{2,p}(\Omega_{\vm})}&\le C(n,d)\|f\|_{W^{n,p}(\Omega_{\vm})}N^{-2(n-2)/d}L^{-2(n-2)/d},\end{align}
		where $C$ is a constant independent of $N$ and $L$.
		
		By combining (iii) and (iv) and setting $K=\lfloor N^{1/d}\rfloor^2\lfloor L^{2/d}\rfloor$, we obtain that for each $\vm\in\{1,2\}^d$, there exists a DSRN $\gamma_{\vm}$ with width $\fO(N\log N)$ and depth $\fO(L\log L)$ such that \begin{align}\|f(\vx)-\gamma_{\vm}(\vx)\|_{W^{2,p}(\Omega_{\vm})}&\le C_3(n,d,a,b)\|f\|_{W^{n,p}(\Omega_{\vm})}N^{-2(n-2)/d}L^{-2(n-2)/d},\end{align}
		where $C_3$ is a constant independent of $N$ and $L$.
		
		\textbf{(v)}: The final step is to combine the sequences $\{s_{\vm}\}_{\vm\in\{1,2\}^d}$ and $\{\gamma_{\vm}\}_{\vm\in\{1,2\}^d}$ to construct a network that can approximate $f$ over the entire space $[a,b]^d$. We define the sequence $\{s_{\vm}\}_{\vm\in\{1,2\}^d}$ because $\gamma_{\vm}$ may not accurately approximate $f$ on $[a,b]^d\backslash\Omega_{\vm}$. The purpose of $s_{\vm}$ is to remove this portion of the domain and to allow other networks to approximate $f$ on $[a,b]^d\backslash\Omega_{\vm}$. 

}

\subsection{Propositions and Lemmas of $\sigma_i$ neural networks}

First, we list a few basic lemmas of $\sigma_i$ neural networks repeatedly applied in our main analysis.
\begin{lemma}\label{sigma2}
The following basic lemmas of $\sigma_2$ neural networks should hold:

(i) $f(x)=x^2$ can be realized exactly by a $\sigma_2$ neural network with one hidden layer and two neurons, and $f(x, y)=x y=\frac{(x+y)^2-(x-y)^2}{4}$ can be realized exactly by a $\sigma_2$ neural network with one hidden layer and four neurons. $x$ can be realized exactly by a $\sigma_2$ neural network with one hidden layer and four neurons or $\sigma_1$ neural network with one hidden layer and two neurons.

(ii) Let \(P(\vx)\) be a polynomial of total degree \(m\) in \(\mathbb{R}^d\). There is a $\sigma_2$-NN of depth \(\lceil\log_2 m\rceil\) and width \(2\,(m+1)^d\) that exactly represents the polynomial \(P(\vx)\).
				
(iii) For any $s>0$, there is a $\sigma_2$-NN with width $2(s+1)$ and depth $\log_2 m$ that exactly represents \( \operatorname{ReLU}^{s+2}\).
			\end{lemma}	

\begin{proof}
(i) is easy to check by $x^2=\sigma_2(x)+\sigma_2(-x)$, $\frac{1}{2}((x+1)^2-x^2-1)=x$ and $x=\sigma_1(x)+\sigma_1(-x)$. (ii) is the results in \cite{he2023expressivity}. For (iii), we know that \begin{equation}
    \operatorname{ReLU}^{s+2}(x)=x^s\cdot\operatorname{ReLU}^2(x).
\end{equation}Then combining (i) and (ii), we have that $ \operatorname{ReLU}^{s+2}$ is a $\sigma_2$ neural network with width $2(s+1)$ and depth $\log_2 m$.
\end{proof}

Define a sequence of subsets of $\Omega$.
\begin{definition}\label{omega}Given $K,d\in\sN^+$, and for any $\vm=(m_1,m_2,\ldots,m_d)\in\{1,2\}^d$, we define $
				\Omega_{\vm}:=\prod_{j=1}^d\Omega_{m_j},
		$  where \begin{align}
			\Omega_1&:=\bigcup_{i=0}^{K-1}\left[a+\frac{i\cdot(b-a)}{K},a+\frac{i\cdot(b-a)}{K}+\frac{3\cdot(b-a)}{4K}\right],\notag\\\Omega_2&:=\bigcup_{i=0}^{K}\left[\frac{i\cdot(b-a)}{K}-\frac{b-a}{2K},a+\frac{i\cdot(b-a)}{K}+\frac{b-a}{4K}\right]\cap [a,b].\end{align}
\end{definition}

   Next, we define a function which will be repeatly used in the proof in this section.
	\begin{definition}
				Define $s(x)$ from $\sR\to [0,1]$ as \begin{equation}
					s(x):=
					\begin{cases}
						
						2x^2,~&x\in\left[0,\frac{1}{2}\right]\\
						-2(x-1)^2+1,~&x\in\left[\frac{1}{2},1\right]\\1,~&x\in \left[1,2\right]  \\-2(x-2)^2+1,~&x\in\left[2,\frac{5}{2}\right]\\
						2(x-3)^2,~&x\in\left[\frac{5}{2},3\right]\\
						0,~&\text{otherwise}.
					\end{cases}
				\end{equation}
	\end{definition}

\begin{definition}\label{sm}
				Given $K\in\sN_+$, we define two functions in $\sR$:\begin{align}
					s_1(x)=\sum_{i=0}^{K}s\left(4Kx+1-4i\right),~s_2(x)=s_1\left(x+\frac{1}{2K}\right).
				\end{align}
				Then for any $\vm=(m_1,m_2,\ldots,m_d)\in\{1,2\}^d$, we define \begin{equation}
					s_{\vm}(\vx):=\prod_{j=1}^d s_{m_j}\left(\frac{x_j-a}{b-a}\right)
				\end{equation} for any $\vx=(x_1,x_2,\ldots,x_d)\in\sR^d$.
\end{definition}

\begin{proposition}\label{smm}
				Given $N,L,d\in\sN_+$ with $\log _2 N>d\log d$, and the setting $K=\lfloor N^{1/d}\rfloor\lfloor L^{2/d}\rfloor$, $\{s_{\vm}(\vx)\}_{\vm\in\{1,2\}^d}$ defined in Definition \ref{sm} satisfies:
				
				(i): $\|s_{\vm}(\vx)\|_{L^\infty([a,b]^d)}\le 1$,~$\|s_{\vm}(\vx)\|_{W^{1,\infty}([a,b]^d)}\le \frac{8K}{b-a}$ and $\|s_{\vm}(\vx)\|_{W^{2,\infty}([a,b]^d)}\le \frac{64K^2}{(b-a)^2}$ for any $\vm\in\{1,2\}^d$.
				
				(ii): $\{s_{\vm}(\vx)\}_{\vm\in\{1,2\}^d}$ is a partition of the unity $[a,b]^d$ with ${\rm supp}~s_{\vm}(\vx)\cap[a,b]^d=\Omega_{\vm}$ defined in Definition \ref{omega}.
				
				(iii): For any $\vm\in\{1,2\}^d$, we have that there is a $\sigma_1$ neural network $\lambda_{\vm,1}(\vx)$ with the width $4N$ and depth $\lceil2\log_2L\rceil$ and a $\sigma_2$ neural network $\lambda_{\vm,2}(\vx)$ with the width $6N$ and depth $\lceil1+\log_2d\rceil$, such as \[\lambda_{\vm,2}\circ\lambda_{\vm,1}(\vx)=\prod_{j=1}^ds_{m_j}\left(\frac{x_j-a}{b-a}\right)=s_{\vm}(\vx),\vx\in[a,b]^d.\]
			\end{proposition}	

   \begin{proof}
				(i) and (ii) are proved by direct calculation. The proof of (iii) follows:
				
				First, we architect $s(x)$ by a $\sigma_2$ neural network. By direct calculation, we notice that \begin{equation}
s(x)=2[x]_{+}^2-4[x-0.5]_{+}^2+2[x-1]_{+}^2-2[x-2]_{+}^2+4[x-2.5]_{+}^2-2[x-3]_{+}^2.
\end{equation} which is a $\sigma_2$ neural network with with width 6 and one hidden layer. The $\widetilde{g}(x)$ defined as \begin{equation}
					\widetilde{g}(x)=\sum_{i=0}^{\lfloor N^{1/d}\rfloor-1}s\left(4Kx-4i-\frac{1}{2}\right)
				\end{equation}is a $\sigma_2$ neural network with the width $12(\lfloor N^{1/d}\rfloor)$ and two hidden layers.

                Next, we construct \( \psi_i \) for \( i = 2, \ldots, \lceil\frac{2}{d}\log_2 L\rceil + 2 \) based on the symmetry and periodicity of \( \widetilde{g}(x) \) due to the symmetry of \( s(x) \). $\psi_2$ is the function with period $\frac{2}{\lfloor N^{1/d}\rfloor\lfloor L^{2/d}\rfloor}$ in $\left[0,\frac{1}{\lfloor L^{2/d}\rfloor}\right]$, and each period is a hat function with gradient 1. $\psi_i$ for $i=3,\ldots,\lceil\frac{2}{d}\log_2 L\rceil+2$ is the hat function in $\left[0,\frac{2^{i-2}}{\lfloor L^{2/d}\rfloor}\right]$. The schematic diagram is in Fig.~\ref{psi4} (The diagram is shown the case for $\lfloor N^{1/d}\rfloor$ is a even integer).

                \begin{figure}[h!]
				\centering
				\includegraphics[scale=0.42]{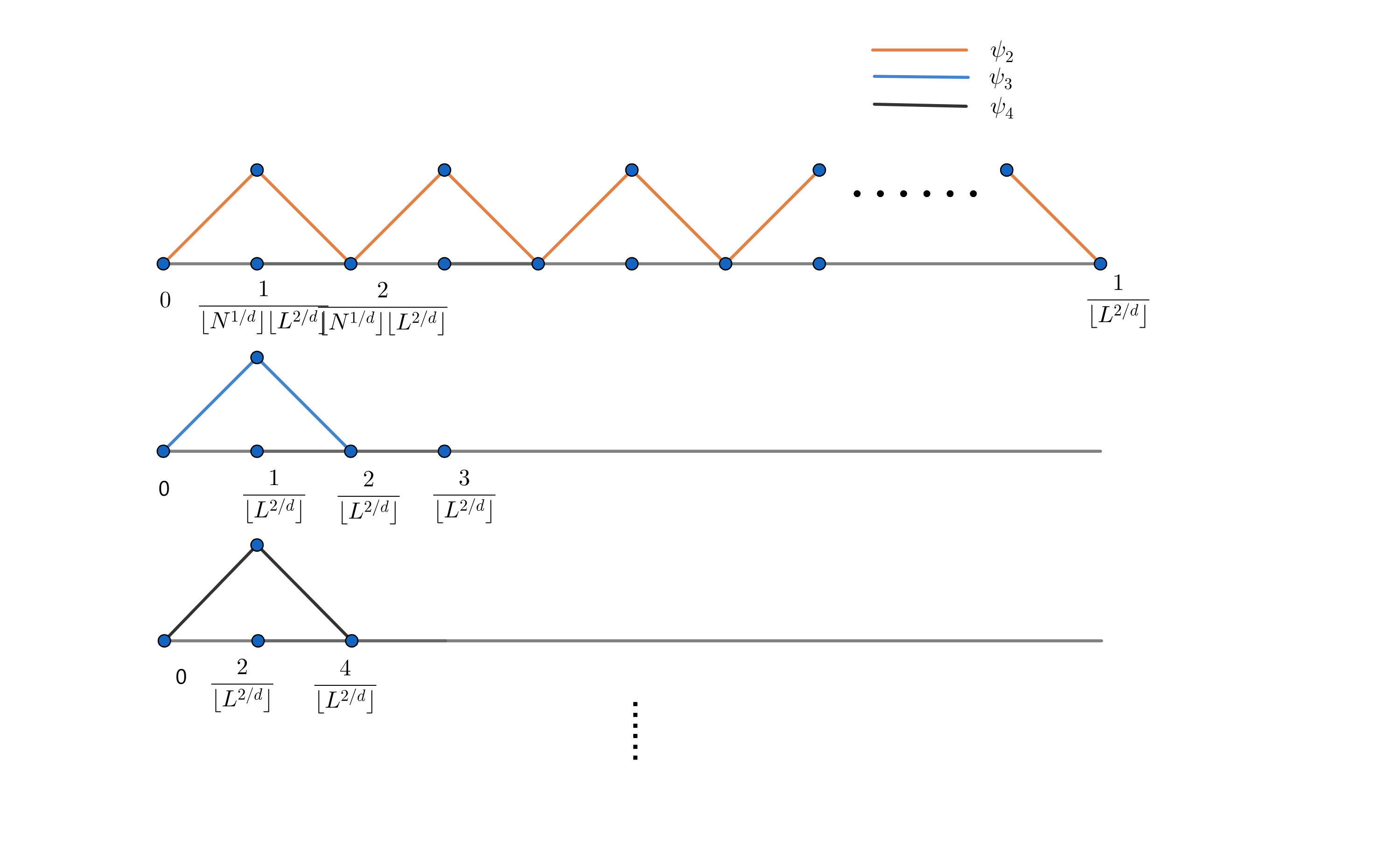}
				\caption{$\psi_i$ for $i=2,\ldots,\lceil\frac{1}{d}\log_2 L\rceil$.}
				\label{psi4}
			\end{figure}

                Note that $\psi_2\circ\psi_3\ldots\circ\psi_{\lceil\frac{2}{d}\log_2 L\rceil+2}(x)$ is the function with period $\frac{2}{\lfloor N^{1/d}\rfloor\lfloor L^{2/d}\rfloor}$ in $[0,1]\subset\left[0,\frac{2^{\lceil\frac{2}{d}\log_2 L\rceil}}{\lfloor L^{2/d}\rfloor}\right]$, and each period is a hat function with gradient 1. Then function \begin{equation}
				\psi(x):=\widetilde{g}(x)\circ \psi_2\circ\psi_3\ldots\circ\psi_{\lceil\frac{2}{d}\log_2 L\rceil+2}(x):=\widetilde{g}\circ\psi_1\label{rep}
			\end{equation} is obtained by repeating reflection $\widetilde{g}(x)$ in $[0,1]$.
			
			$\psi_2$ is a network with width $4\lfloor N^{1/d}\rfloor$ and one hidden layer. Rest $\psi_i$ is a $\sigma_1$-NN with width 3 and one hidden layer. Hence $\psi_1(x)$ is a $\sigma_1$-NN with width $4\lfloor N^{1/d}\rfloor$ and depth $\lceil\frac{2}{d}\log_2 L\rceil$ and $s_1(x)=\psi\left(x+\frac{1}{8K}\right)$ and $s_2(x)=\psi\left(x+\frac{5}{8K}\right)$.

            Based on (ii) in Lemma \ref{sigma2}, we have that there is a $\sigma_1$ neural network $\lambda_{\vm,1}(\vx)$ with the width $4N$ and depth $\lceil2\log_2L\rceil$ and a $\sigma_2$ neural network $\lambda_{\vm,2}(\vx)$ with the width $6N$ and depth $\lceil1+\log_2d\rceil$, such as \[\lambda_{\vm,2}\circ\lambda_{\vm,1}(\vx)=\prod_{j=1}^ds_{m_j}\left(\frac{x_j-a}{b-a}\right)=s_{\vm}(\vx),\vx\in[a,b]^d,\]where the order of width is based on $dN^{1/d}\le N$ for $N>d$.
            \end{proof}

Next, we present some propositions for ReLU DNNs.
		
		\begin{proposition}[\cite{lu2021deep}]\label{step}
			Given any $N,L\in\sN_+$ and $\delta\in\Big(0,\frac{1}{3K}\Big]$ for $K=\lfloor N^{1/d}\rfloor^2\lfloor L^{2/d}\rfloor$, there exists a $\sigma_1$-NN $\phi$ with the width $4N+5$ and depth $4L+4$ such that 
			
			\[\phi(x)=k,x\in\left[a+\frac{k\cdot(b-a)}{K},a+\frac{(k+1)\cdot(b-a)}{K}-\delta\cdot 1_{k< a+(K-1)\cdot(b-a)}\right], \]for \(k=0,1,\dots,K-1.\)
		\end{proposition}
        \begin{remark}
            In \cite{lu2021deep}, this proposition is proved on the \([0,1]\) domain. Here, we extend it to the \([a,b]\) domain, which is straightforward by applying a rescaling and transformation.
        \end{remark}

		\begin{proposition}[\cite{lu2021deep}]\label{point}
			Given any $N,L,s\in\sN_+$ and $\xi_i\in[0,1]$ for $i=0,1,\ldots N^2L^2-1$, there exists a $\sigma_1$-NN $\phi$ with the width $16s(N+1)\log_2(8N)$ and depth $(5L+2)\log_2(4L)$ such that
			
			1. $|\phi(i)-\xi_i|\le N^{-2s}L^{-2s}$ for $i=0,1,\ldots N^2L^2-1$.
			
			2. $0\le \phi(x)\le 1$, $x\in\sR$.
		\end{proposition}

\subsection{Proof of Theorem \ref{main2}}  
To derive the approximation result, it can be broken down into three sequential steps. Without sacrificing generality, we focus on approximating functions in $W^{n,p}$ using the norm measured in $W^{m,p}$.

\subsubsection{Proof of Theorem \ref{main2}}
Define subsets of $\Omega_{\vm}$ for simplicity notations.
		
For any $\vm\in\{1,2\}^d$, we define \begin{equation}
			\Omega_{\vm,\vi}:=[a,b]^d\cap\prod_{j=1}^d\left[a+\frac{2i_j-1_{m_j< 2}}{2K}\cdot(b-a),a+\frac{3+4i_j-2\cdot1_{m_j< 2}}{4K}\cdot(b-a)\right]
		\end{equation}$\vi=(i_1,i_2,\ldots,i_d)\in\{0,1\ldots,K\}^d$, and it is easy to check $\bigcup_{\vi\in\{0,1\ldots,K\}^d}\Omega_{\vm,\vi}=\Omega_{\vm}$. 
  
  \begin{proposition}[{\cite[Theorem 8]{yang2023nearly}}]\label{fN}
			Let $K\in\sN_+$, $p\in[1,\infty]$ and $n\ge 2$. Then for any $f\in W^{n,p}([a,b]^d)$ and $\vm\in\{1,2\}^d$, there exist piece-wise polynomials functions $f_{K,\vm}=\sum_{|\valpha|\le n-1}g_{f,\valpha,\vm}(\vx)\vx^{\valpha}$ on $\Omega_{\vm}$ (Definition \ref{omega}) with the following properties:
			\begin{align}\|f-f_{K,\vm}\|_{W^{2,p}(\Omega_{\vm})}&\le C_1(n,d,|b-a|)\|f\|_{W^{n,p}(\Omega_{\vm})}K^{-(n-2)},\notag\\
					\|f-f_{K,\vm}\|_{W^{1,p}(\Omega_{\vm})}&\le C_1(n,d,|b-a|)\|f\|_{W^{n,p}(\Omega_{\vm})}K^{-(n-1)},\notag\\\|f-f_{K,\vm}\|_{L^{p}(\Omega_{\vm})}&\le C_1(n,d,|b-a|)\|f\|_{W^{n,p}(\Omega_{\vm})}K^{-n}.  
				\end{align}Furthermore, $g_{f,\valpha,\vm}(\vx):\Omega_{\vm}\to\sR$ is a constant function on each $\Omega_{\vm,\vi}$ for $\vi\in\{0\ldots,K\}^d,$ and \begin{equation}
				|g_{f,\valpha,\vm}(\vx)|\le C_2(n,d,|a|,|b|)\|f\|_{W^{n-1,p}(\Omega_{\vm})}(4K)^{\frac{d}{p}}
			\end{equation} for all $\vx\in\Omega_{\vm}$, where $C_1$ and $C_2$ are constants independent of $K$.
		\end{proposition}
		
		This proof is similar to that of \cite[Lemma C.4]{guhring2020error} and \cite[Theorem 8]{yang2023nearly}, so we omit detailed proof in this paper.

      \begin{remark}\label{C}
The constant \(C_1\) in Proposition~\ref{fN} depends linearly on \(\max\{|b - a|^{n - i}\}_{i = 0}^2\). This dependence arises from the Bramble–Hilbert Lemma \cite[Lemma 4.3.8]{brenner2008mathematical}. As shown in the proof of \cite[Theorem 8]{yang2023nearly} and \cite[Lemma 4.3.8]{brenner2008mathematical}, the approximation constant in the Bramble–Hilbert Lemma scales with \((\operatorname{diam}(\Omega))^s\), where \(\Omega\) is the approximation domain and \(s\) is the difference between the regularity of the function and the order of the norm in which the error is measured. The constant \(C_2\) depends on the maximal distance between the approximation domain and the origin. Since the target function is approximated locally by polynomials, similar to a Taylor expansion, the coefficients of these polynomials naturally depend on the location within the domain. In our setting, \(g_{f,\valpha,\vm}\) is a piecewise constant function that represents the local coefficients of the polynomial approximation. The bound on \(g_{f,\valpha,\vm}\), and hence the constant \(C_2\), depends on \(\max\{|a|, |b|\}\), and the degree of this dependence is determined by the order of the polynomial. A more detailed analysis can be found in \cite[Lemma B.9]{guhring2020error}.
\end{remark}

			\begin{proof}[Proof of Theorem \ref{main2}] 
				Without loss of the generalization, we consider the case for $\vm_*=(1,1,\ldots,1)$.
				Due to Proposition \ref{fN} and setting $K=\lfloor N^{1/d}\rfloor^2\lfloor L^{2/d}\rfloor$, we have \begin{align}&\|f-f_{K,\vm_*}\|_{W^{2,p}(\Omega_{\vm_*})}\le C_1(n,d,|b-a|)\|f(\vx)\|_{W^{n,p}(\Omega_{\vm_*})}N^{-2(n-2)/d}L^{-2(n-2)/d}\notag\\&\|f-f_{K,\vm_*}\|_{W^{1,p}(\Omega_{\vm_*})}\le C_1(n,d,|b-a|)\|f(\vx)\|_{W^{n,p}(\Omega_{\vm_*})}N^{-2(n-1)/d}L^{-2(n-1)/d}\notag\\&\|f-f_{K,\vm_*}\|_{L^{p}(\Omega_{\vm_*})}\le C_1(n,d,|b-a|)\|f(\vx)\|_{W^{n,p}(\Omega_{\vm_*})}N^{-2n/d}L^{-2n/d},\label{thm3}
				\end{align} where $f_{K,\vm_*}=\sum_{|\valpha|\le n-1}g_{f,\valpha,\vm_*}(\vx)\vx^{\valpha}$ for $x\in \Omega_{\vm_*}$. Note that $g_{f,\valpha,\vm_*}(\vx)$ is a constant function for $\vx\in\prod_{j=1}^d\left[a+\frac{i_j}{K}\cdot(b-a),a+\frac{3+4i_j}{4K}\cdot(b-a)\right]$ and $\vi=(i_1,\ldots,i_d)\in\{0,1,\ldots,K-1\}^d$. The remaining part is to approximate $f_{K,\vm_*}$ by neural networks.
				
				The way to approximate $g_{f,\valpha,\vm_*}(\vx)$ is similar with \cite[Theorem 3.1]{hon2022simultaneous}. First of all, due to Proposition \ref{step}, there is a neural network $\phi_1(x)$ with the width $4N+5$ and depth $4L+4$ such that\begin{equation}
				\phi(x)=k,x\in\left[a+\frac{k\cdot(b-a)}{K},a+\frac{(k+1)\cdot(b-a)}{K}-\frac{b-a}{4K}\right], ~k=0,1,\ldots,K-1.
			\end{equation} Note that we choose $\delta=\frac{1}{4K}\le \frac{1}{3K}$ in Proposition \ref{step}. Then define \[\vphi_2(\vx)=\left[\frac{\phi_1(x_1)}{K},\frac{\phi_1(x_2)}{K},\ldots,\frac{\phi_1(x_d)}{K}\right]^\T.\] For each $p=0,1,\ldots,K^d-1$, there is a bijection\[\veta(p)=[\eta_1,\eta_2,\ldots,\eta_d]\in \{0,1,\ldots,K-1\}^d\] such that $\sum_{j=1}^d\eta_jK^{j-1}=p$. Then define \[\xi_{\valpha,p}=\frac{g_{f,\valpha,\vm_*}\left(\frac{\veta(p)}{K}\right)+C_2(n,d,|a|,|b|)\|f\|_{W^{n-1,p}(\Omega_{\vm_*})}(4K)^{\frac{d}{p}}}{2C_2(n,d,|a|,|b|)\|f\|_{W^{n-1,p}(\Omega_{\vm_*})}(4K)^{\frac{d}{p}}}\in[0,1],\]where $C_2(n,d,|a|,|b|)\|f\|_{W^{n-1,p}(\Omega_{\vm_*})}(4K)^{\frac{d}{p}}$ is the bounded of $g_{f,\valpha,\vm_*}$ and $C_2$ defined in Proposition \ref{fN}.
			Therefore, based on Proposition \ref{point}, there is a neural network $\tilde{\phi}_{\valpha}(x)$ with the width \[16\left(n+\frac{d}{p}\right)(N+1)\log_2(8N)\] and depth $(5L+2)\log_2(4L)$ such that $|\tilde{\phi}_{\valpha}(p)-\xi_{\valpha,p}|\le N^{-2\left(n-\frac{d}{p}\right)}L^{-2\left(n-\frac{d}{p}\right)}$ for $p=0,1,\ldots K^d-1$. Denote \[\phi_{\valpha}(\vx)=C_2(n,d,|a|,|b|)(4K)^{\frac{d}{p}}\|f\|_{W^{n-1,p}(\Omega_{\vm_*})}\left[2\tilde{\phi}_{\valpha}\left(\sum_{j=1}^d\eta_jK^j\right)-1\right]\] and obtain that\begin{align}&\left|\phi_{\valpha}\left(\frac{\veta(p)}{K}\right)-g_{f,\valpha,\vm_*}\left(\frac{\veta(p)}{K}\right)\right|\notag\\=&2C_2(n,d,|a|,|b|)\|f\|_{W^{n-1,p}(\Omega_{\vm_*})}(4K)^{\frac{d}{p}}|\tilde{\phi}_{\valpha}(p)-\xi_{\valpha,p}|\notag\\\le & 4^{\frac{2d+p}{2p}}\|f\|_{W^{n-1,p}(\Omega_{\vm_*})}C_2(n,d,|a|,|b|)(NL)^{-2n}.\notag\end{align} Then we obtain that \begin{align}
				\|\phi_{\valpha}\left(\vphi_2(\vx)\right)-g_{f,\valpha,\vm_*}\left(\vx\right)\|_{W^{1,p}(\Omega_{\vm_*})}=&\|\phi_{\valpha}\left(\vphi_2(\vx)\right)-g_{f,\valpha,\vm_*}\left(\vx\right)\|_{L^{p}(\Omega_{\vm_*})}\notag\\\le&4^{\frac{2d+p}{2p}}C_2(n,d,|a|,|b|)\|f\|_{W^{n-1,p}(\Omega_{\vm_*})}(NL)^{-2n}
			\end{align}which is due to $\phi_{\valpha}\left(\vphi_2(\vx)\right)-g_{f,\valpha,\vm_*}\left(\vx\right)$ is a step function, and the first order weak derivative is $0$ in $\Omega_{\vm_*}$, and we have that \begin{align}
					\|\phi_{\valpha}\left(\vphi_2\right)-g_{f,\valpha,\vm_*}\left(\vx\right)\|_{W^{2,p}(\Omega_{\vm_*})}=&\|\phi_{\valpha}\left(\vphi_2\right)-g_{f,\valpha,\vm_*}\left(\vx\right)\|_{W^{1,p}(\Omega_{\vm_*})}\notag\\=&\|\phi_{\valpha}\left(\vphi_2\right)-g_{f,\valpha,\vm_*}\left(\vx\right)\|_{L^{p}(\Omega_{\vm_*})}\notag\\\le&4^{\frac{2d+p}{2p}}C_2(n,d,|a|,|b|)\|f\|_{W^{n-1,p}(\Omega_{\vm_*})}(NL)^{-2n}.\label{piece}
				\end{align}
				
				Due to (v) in Lemma \ref{sigma2}, there is a $\sigma_2$ neural network $\phi_{5,\valpha}(\vx)$ with the width $N$ due to $N>n^d$ and depth $\left\lceil\log _2 L\right\rceil$ due to $L>n$ such that \begin{equation}
					\phi_{5,\valpha}(\vx)=\vx^{\valpha},~\vx\in\sR^d.
				\end{equation} Due to (iv) in Lemma \ref{sigma2}, there is a $\sigma_2$ neural network $\phi_{6}(\vx)$ with the width $4$ and depth $1$ such that \begin{equation}
					\phi_{6}(x,y)=xy,~x,y\in\sR.
				\end{equation}
				
				Now we define the neural network $\gamma_{\vm_*}(\vx)$ to approximate $f_{K,{\vm_*}}(\vx)$ in $\Omega_{\vm_*}$:\begin{equation}
					\gamma_{\vm_*}(\vx)=\sum_{|\valpha|\le n-1}\phi_6\left[\phi_{\valpha}(\vphi_2(\vx)),\phi_{5,\valpha}(\vx)\right].\label{2nnphi}
				\end{equation} $\gamma_{\vm_*}$ is a $\sigma_2$-NN with depth $\log_29(L+1)$ and the width $28n^{d+1}(N+d)\log_2(8N)$.
    
            The remaining question is to find the error $\fE$:
				\begin{align}
					\widetilde{\fE}:=&\left\|\sum_{|\valpha|\le n-1}\phi_6\left[\phi_{\valpha}(\vphi_2(\vx)),\phi_{5,\valpha}(\vx)\right]-f_{K,\vm_*}(\vx)\right\|_{W^{2,p}(\Omega_{\vm_*})}\notag\\\le&\sum_{|\valpha|\le n-1}\left\|\phi_6\left[\phi_{\valpha}(\vphi_2(\vx)),\phi_{5,\valpha}(\vx)\right]-g_{f,\valpha,\vm_*}(\vx)\vx^{\valpha}\right\|_{W^{2,p}(\Omega_{\vm_*})}\notag\\=&\sum_{|\valpha|\le n-1}\left\|\phi_{\valpha}(\vphi_2(\vx))\vx^{\valpha}-g_{f,\valpha,\vm_*}(\vx)\vx^{\valpha}\right\|_{W^{2,p}(\Omega_{\vm_*})}\notag\\\le&n^2\sum_{|\valpha|\le n-1}\left\|\phi_{\valpha}(\vphi_2(\vx))-g_{f,\valpha,\vm_*}(\vx)\right\|_{W^{2,p}(\Omega_{\vm_*})}\notag\\\le &n^{d+2}4^{\frac{2d+p}{2p}}C_2(n,d,|a|,|b|)\|f\|_{W^{n-1,p}(\Omega_{\vm_*})}(NL)^{-2n}.\label{thm4}
				\end{align}

				Combining Eqs.~(\ref{thm3}) and (\ref{thm4}), we have that there is a DSRN $\gamma_{\vm_*}$ belonging $\fN_{2,9(L+1)\log_2L}$ with the width $28n^{d+1}(N+d)\log_2(8N)$ such that\begin{align}&\|f(\vx)-\gamma_{\vm_*}(\vx)\|_{W^{2,p}(\Omega_{\vm_*})}\le C_3(n,d,a,b)\|f\|_{W^{n,p}(\Omega_{\vm_*})}N^{-2(n-2)/d}L^{-2(n-2)/d}\notag\\&\|f(\vx)-\gamma_{\vm_*}(\vx)\|_{W^{1,p}(\Omega_{\vm_*})}\le C_3(n,d,a,b)\|f\|_{W^{n,p}(\Omega_{\vm_*})}N^{-2(n-1)/d}L^{-2(n-1)/d}\notag\\&\|f(\vx)-\gamma_{\vm_*}(\vx)\|_{L^{p}(\Omega_{\vm_*})}\le C_3(n,d,a,b)\|f\|_{W^{n,p}(\Omega_{\vm_*})}N^{-2n/d}L^{-2n/d},\label{app2}\end{align}
				where $C_3=C_1+4^{\frac{2d+p}{2p}}C_2$ is the constant independent of $N,L$. The final inequality arises from the fact that $\frac{2n}{d}\le 2n$ holds for values of $d$ and $n$ greater than or equal to 1.
				
				Similarly, we can construct a $\gamma_{\vm}$ achieving the same order of Eq.~({\ref{app2}}).

                Denote a sequence of the neural network to achieve performance in Proposition \ref{smm} as $\{s_{\vm}(\vx):=\lambda_{\vm,2}\circ\lambda_{\vm,1})(\vx)\}_{\vm\in\{1,2\}^d}$.
				
				Now we define the piece-wise constant function $\phi_{\valpha}(\vphi_2(\vx))$ (\ref{piece}) defined on $\Omega_{\vm}$ as $p_{\vm}$, then we have \begin{align}
					\gamma(\vx)&=\sum_{\vm\in\{1,2\}^d}s_{\vm}(\vx)\cdot\gamma_{\vm}(\vx)=\sum_{\vm\in\{1,2\}^d}s_{\vm}(\vx)\cdot\sum_{|\valpha|\le n-1}p_{\vm}\cdot\vx^{\valpha}\notag\\&=\sum_{|\valpha|\le n-1}\sum_{\vm\in\{1,2\}^d}s_{\vm}(\vx)\cdot p_{\vm}\cdot\vx^{\valpha}\notag\\&=\sum_{|\valpha|\le n-1}\sum_{\vm\in\{1,2\}^d}(\lambda_{\vm,2}\circ\lambda_{\vm,1})(\vx)\cdot p_{\vm}\cdot\vx^{\valpha}\notag\\&=\sum_{|\valpha|\le n-1}\sum_{\vm\in\{1,2\}^d}\phi_6(\phi_6((\lambda_{\vm,2},\textbf{Id})\circ(\lambda_{\vm,1},p_{\vm})),\vx^{\valpha}),\label{4nnphi}
				\end{align} where we denote $(f_1,f_2)\circ(g_1,g_2)$ as $(f_1\circ g_1,f_2\circ g_2)$. It is easy to check that $\gamma$ is a DSRN in $\fN_{1,9(L+1)\log_24L}$ based on Fig.~\ref{gammaaa}. Furthermore, the width of $\gamma_{\vm_*}$ is $2^{d+6}n^{d+1}(N+d)\log_2(8N)$.

    \begin{figure}[h!]
				\centering
				\includegraphics[scale=0.07]{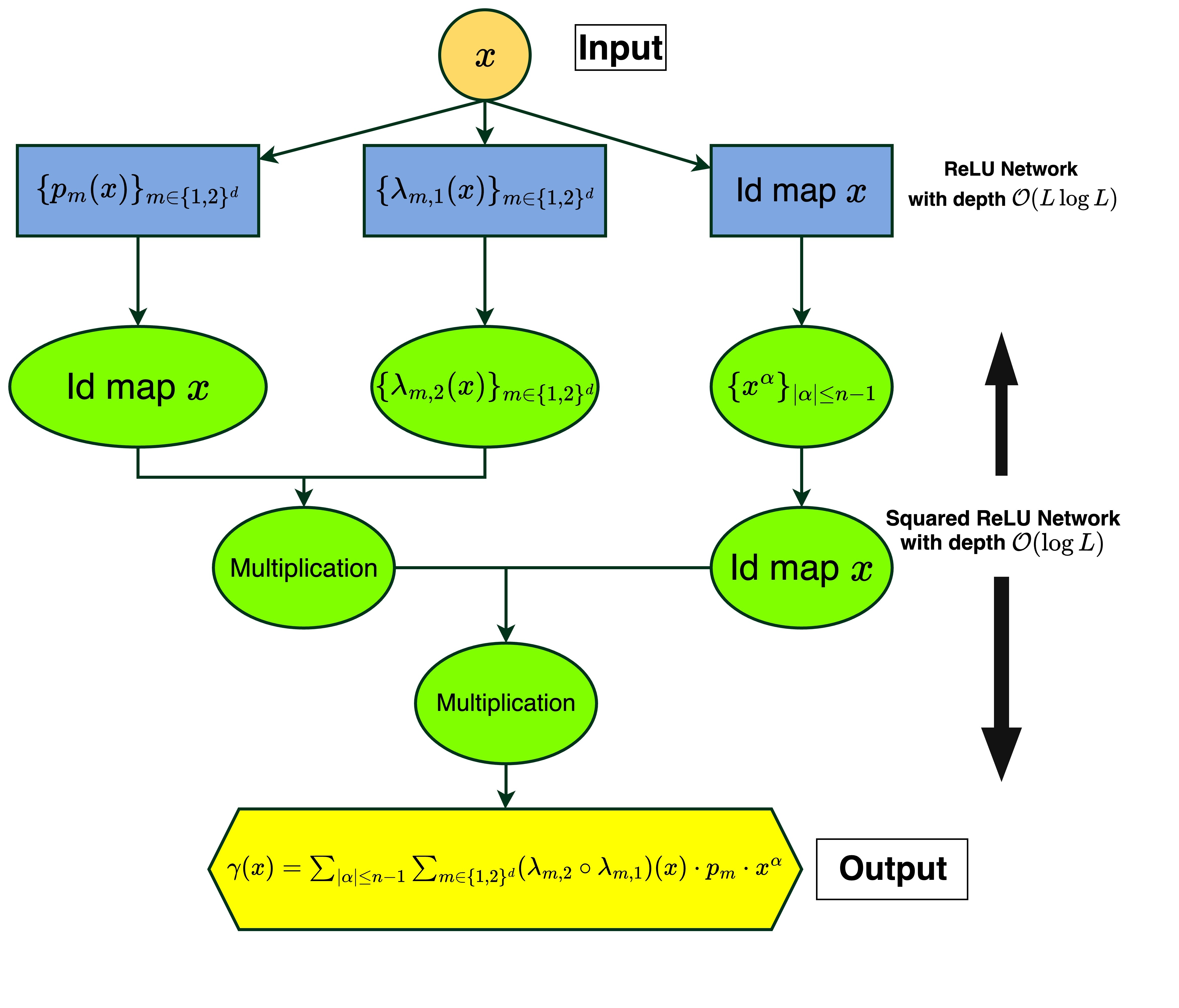}
				\caption{The structure of $\gamma(\vx)$.}
				\label{gammaaa}
			\end{figure}
				
				Note that \begin{align}
					\widetilde{\fR}:=&\|f(\vx)-\gamma(\vx)\|_{W^{2,p}([a,b]^d)}\le\sum_{\vm\in\{1,2\}^d}\left\| s_{\vm}(x)\cdot f(\vx)-s_{\vm}(\vx)\gamma_{\vm}(\vx)\right\|_{W^{2,p}([a,b]^d)}\notag\\= &\sum_{\vm\in\{1,2\}^d}\left\| s_{\vm}(x)\cdot f(\vx)-s_{\vm}(\vx)\gamma_{\vm}(\vx)\right\|_{W^{2,p}(\Omega_{\vm})},\notag
				\end{align}
    where the last equality is due to ${\rm supp}~s_{\vm}(\vx)\cap[a,b]^d=\Omega_{\vm}$.
				
				Then due to chain rule, for each $\vm\in\{1,2\}^d$, we have\begin{align}
&\bigl\|\,s_{\vm}f-s_{\vm}\gamma_{\vm}\bigr\|_{W^{2,p}(\Omega_{\vm})}\notag\\
\le{}&\|s_{\vm}\|_{W^{2,\infty}(\Omega_{\vm})}\|f-\gamma_{\vm}\|_{L^{p}(\Omega_{\vm})}
      +2\|s_{\vm}\|_{W^{1,\infty}(\Omega_{\vm})}\|f-\gamma_{\vm}\|_{W^{1,p}(\Omega_{\vm})}\notag\\
     &+\|s_{\vm}\|_{L^{\infty}(\Omega_{\vm})}\|f-\gamma_{\vm}\|_{W^{2,p}(\Omega_{\vm})}
      +\|s_{\vm}\|_{W^{1,\infty}(\Omega_{\vm})}\|f-\gamma_{\vm}\|_{L^{p}(\Omega_{\vm})}\notag\\
     &+\|s_{\vm}\|_{L^{\infty}(\Omega_{\vm})}\|f-\gamma_{\vm}\|_{W^{1,p}(\Omega_{\vm})}
      +\|s_{\vm}\|_{L^{\infty}(\Omega_{\vm})}\|f-\gamma_{\vm}\|_{L^{p}(\Omega_{\vm})}\notag\\
\le{}&91\,C_3(n,d,a,b)\,\|f\|_{W^{n,p}([a,b]^{d})}(NL)^{-2(n-2)/d}.\notag
\end{align}

				Hence \[\widetilde{\fR}\le 2^{d+7}C_3(n,d,a,b)\|f\|_{W^{n,p}([a,b]^d)}N^{-2(n-2)/d}L^{-2(n-2)/d}.\]
				
			\end{proof}

   Overall, we have shown that DSRNs can approximate functions in Sobolev spaces with the \(W^{2,p}\) norm. Our method extends directly to \(W^{m,p}\) for any integer \(m \ge 1\). The key modification in the proof is the construction of a \(C^m\) partition of unity \(\{s_{\vm}(\vx)\}_{\vm \in \{1,2\}^d}\), which we achieve by designing each \(s_{\vm}(\vx)\) as a piecewise polynomial of sufficiently high degree. In particular, to build a \(C^3\) partition of unity, one may take the one-dimensional “bump” function
\[
s(x) \;=\;
\begin{cases}
0, & x \le 0,\\
6x^5 \;-\; 15x^4 \;+\; 10x^3, & 0 < x < 1,\\
1, & 1 \le x \le 2,\\
6(3 - x)^5 \;-\; 15(3 - x)^4 \;+\; 10(3 - x)^3, & 2 < x < 3,\\
0, & x \ge 3.
\end{cases}
\]
This \(C^3\) function admits an exact representation by \(\mathrm{ReLU}^5\):
\[
s(x)
\;=\;
\sum_{k=0}^5 (-1)^k \binom{5}{k}\,
\operatorname{ReLU}^5\Bigl(\tfrac{x}{3} + \tfrac{2}{3} - k\Bigr).
\]
By Lemma \ref{sigma2}, each \(\mathrm{ReLU}^5\) can itself be realized by a \(\sigma_2\)-NN. Consequently, we obtain a DSRN that approximates the target function in the \(W^{3,p}\) norm. Repeating this construction inductively yields Corollary \ref{main4}.

   \begin{corollary}\label{main4}
			For any $f\in W^{n,p}([a,b]^d)$, $m\in\sN$ with $m\ge 2$ and $1\le p\le +\infty$, any $N, L\in\sN_+$ with $N \log_2L+2^{\left\lfloor\log _2 N\right\rfloor} \geq \max\{d,n\}$ and $L\ge N,m< n$, there is a DSRN $\gamma(\vx)$ in $\fN_{\eta_1,\eta_2L\log_2L}$ with the width $\eta_3N\log_2N$ such that\[\|f(\vx)-\gamma(\vx)\|_{W^{m,p}([a,b]^d)}\le C_4(n,d,a,b)\|f\|_{W^{n,p}([a,b]^d)}N^{-2(n-m)/d}L^{-2(n-m)/d},\]where $\eta_i, C_{11}$ are the constants independent of $N,L$.
\end{corollary}

{In the above proof, we remove the non-smooth points using a smooth partition of unity. In numerical practice, even if a DSRN with a random initial condition does not contain second-order derivatives across the entire domain, it does have second-order derivatives at most points, meaning that the non-smooth points are of zero measure. If we unfortunately sample points at the non-smooth locations, we can resample in the next iteration to avoid these points. In Sec.~\ref{num}, we demonstrate that DSRN works effectively in practice.
}

\section{Optimality of Theorem \ref{main2}}
In this section, our aim is to demonstrate the optimality of DSRN approximations in Sobolev spaces measured by norms in $W^{m,p}\left([a,b]^d\right)$ for $m\geq 2$. Without loss of generality, we prove the optimality for $m=2$.

\begin{theorem}\label{Optimality}
			Given any $\rho, D_{1}, D_{2}, D_{3}, D_4, J_0>0$ and $n,d\in\sN^+$,  there exist $N,L\in\sN$ with $NL\ge J_0$ and $f$, satisfying for any $\phi\in\fN_{D_1,D_{2} L\log L}$ with the width smaller than $D_{3} N\log N$, we have \begin{equation}
				|\phi-f|_{W^{2,p}([a,b]^d)}> D_{4}\|f\|_{W^{n,p}([a,b]^d)}L^{-2(n-2)/d-\rho}N^{-2(n-2)/d-\rho}.
			\end{equation}
		\end{theorem}

%  \subsection{Vapnik--Chervonenkis dimension of second order derivative of DSRNs}

In order to prove the optimality of Theorem \ref{main2}, we need to bound the VC-dimension of second order derivative of deep neural networks (DNNs).

\begin{definition}[VC-dimension \cite{abu1989vapnik}]
	Let $H$ denote a class of functions from $\fX$ to $\{0,1\}$. For any non-negative integer $m$, define the growth function of $H$ as \[\Pi_H(m):=\max_{x_1,x_2,\ldots,x_m\in \fX}\left|\{\left(h(x_1),h(x_2),\ldots,h(x_m)\right): h\in H \}\right|.\] The Vapnik--Chervonenkis dimension (VC-dimension) of $H$, denoted by $\operatorname{VCdim}(H)$, is the largest $m$ such that $\Pi_H(m)=2^m$. For a class $\fG$ of real-valued functions, define $\operatorname{VCdim}(\fG):=\operatorname{VCdim}(\sgn(\fG))$, where $\sgn(\fG):=\{\sgn(f):f\in\fG\}$ and $\sgn(x)=1[x>0]$.
	
\end{definition}

 \begin{theorem}\label{vcdim}
		For any $N,L,d,C\in\sN_+$, there exists a constant $\bar{C}$ independent of $N,L$ such that	\begin{equation}
		\operatorname{VCdim}(D^2\Phi)\le \bar{C} N^2L^2\log_2 L\log_2 N,\label{bound}
       \end{equation}for
       \begin{align}
	D^2\Phi:=\left\{\lambda(\vx)=D^{\valpha}\phi:\phi\in\Phi,~|\valpha|=2\right\},\end{align}where $ \Phi:=\left\{\phi:\phi\text{ is a DSRN in $\fN_{C,L}$ with width$\le N$},\phi\in W^2([a,b]^d)\right\}$, and \[D^{\valpha}=D^{\alpha_1}_1D^{\alpha_2}_2\ldots D^{\alpha_d}_d\] as the partial derivative where $\valpha=[\alpha_{1},\alpha_{2},\ldots,\alpha_d]^T$ and $D_i$ is the weak derivative in the $i$-th variable.
	\end{theorem}

 Prior to proving Theorem \ref{vcdim}, we introduce two lemmas.

 \begin{lemma}[\cite{bartlett2019nearly},\cite{anthony1999neural}]\label{bounded}
			Suppose $W\le M$ and let $P_1,\ldots,P_M$ be polynomials of degree at most $D$ in $W$ variables. Define \[P:=\left|\{\left(\sgn(P_1(\va)),\ldots,\sgn(P_M(\va))\right):\va\in\sR^W\}\right|,\] then we have $P\le 2(2eMD/W)^W$.
		\end{lemma}

   \begin{lemma}[\cite{bartlett2019nearly}]\label{inequality}
			Suppose that $2^m\le 2^t(mr/w)^w$ for some $r\ge 16$ and $m\ge w\ge t\ge0$. Then, $m\le t+w\log_2(2r\log_2r)$.
		\end{lemma}

  \begin{proof}[Proof of Theorem \ref{Optimality}]

The difficulty of the estimate raises from the complexity of the second order derivatives of DNNs: \[\phi=\vW_{L+1}\sigma(\vW_{L}\sigma(\ldots\sigma(\vW_1\vx+\vb_1)\ldots)+\vb_L)+b_{L+1},\]where $\sigma$ can be either the ReLU or the ReLU square. Then the first order derivative can be read as \begin{align}
	\psi(\vx)=D_i \phi(\vx)=&\vW_{L+1}\sigma'(\vW_{L}\sigma(\ldots\sigma(\vW_1\vx+\vb_1)\ldots)+\vb_L)\notag\\&\cdot \vW_{L}\sigma'(\ldots\sigma(\vW_1\vx+\vb_1)\ldots)\ldots\vW_2\sigma'(\vW_1\vx+\vb_1)(\vW_1)_i,
\end{align}where $\vW_i\in\sR^{N_i\times N_{i-1}}$ ($(\vW)_i$ is $i$-th column of $\vW$) and $\vb_i\in\sR^{N_i}$ are the weight matrix and the bias vector in the $i$-th linear transform in $\phi$, and $\sigma'(\vx)=\diag(\sigma'(x_i))$. Then we have \begin{align}
\lambda(\vx)=D_jD_i\phi(\vx)=\sum_{a=1}^{L+1} \lambda_a(\vx).\notag
\end{align} For \begin{align}
&\lambda_a(\vx)=\vW_{L+1}\sigma'(\vW_{L}\sigma(\ldots\sigma(\vW_1\vx+\vb_1)\ldots)+\vb_L)\cdot  \vW_{L}\sigma'(\ldots\sigma(\vW_1\vx+\vb_1)\ldots)\ldots\notag\\&\cdot[\vW_{a}\sigma''(\ldots\sigma(\vW_1\vx+\vb_1)\ldots)\cdot \vW_{a-1}\sigma'(\ldots\sigma(\vW_1\vx+\vb_1)\ldots)\vW_2\sigma'(\vW_1\vx+\vb_1)(\vW_1)_j]\notag\\&\cdot\ldots\vW_{a-1}\sigma'(\ldots\sigma(\vW_1\vx+\vb_1)\ldots)\vW_2\sigma'(\vW_1\vx+\vb_1)(\vW_1)_i,\notag
\end{align}where $\sigma''(\vx)=\diag(\sigma''(x_i))$ is a three-order tensor. The multiplication at here is defined as follows: Let a third-order tensor $\mathcal{T} \in \mathbb{R}^{I \times J \times K}$ and three matrices
$A \in \mathbb{R}^{P \times I}$, 
$B \in \mathbb{R}^{Q \times J}$, 
and 
$C \in \mathbb{R}^{R \times K}$.
The multiplication of the tensor with these matrices is defined as
\[
\mathcal{T}' = \mathcal{T} \times_1 A \times_2 B \times_3 C:=(\mathcal{T}'_{pqr})_{p,q,r} \;\;\in \mathbb{R}^{P \times Q \times R},
\]
where
\[
\mathcal{T}'_{pqr} = \sum_{i=1}^{I} \sum_{j=1}^{J} \sum_{k=1}^{K}
A_{pi}\, B_{qj}\, C_{rk}\, \mathcal{T}_{ijk}.
\]
This operation applies a linear transformation to each mode of the tensor $\mathcal{T}$.

 \begin{figure}[h!]
				\centering
				\includegraphics[scale=0.077]{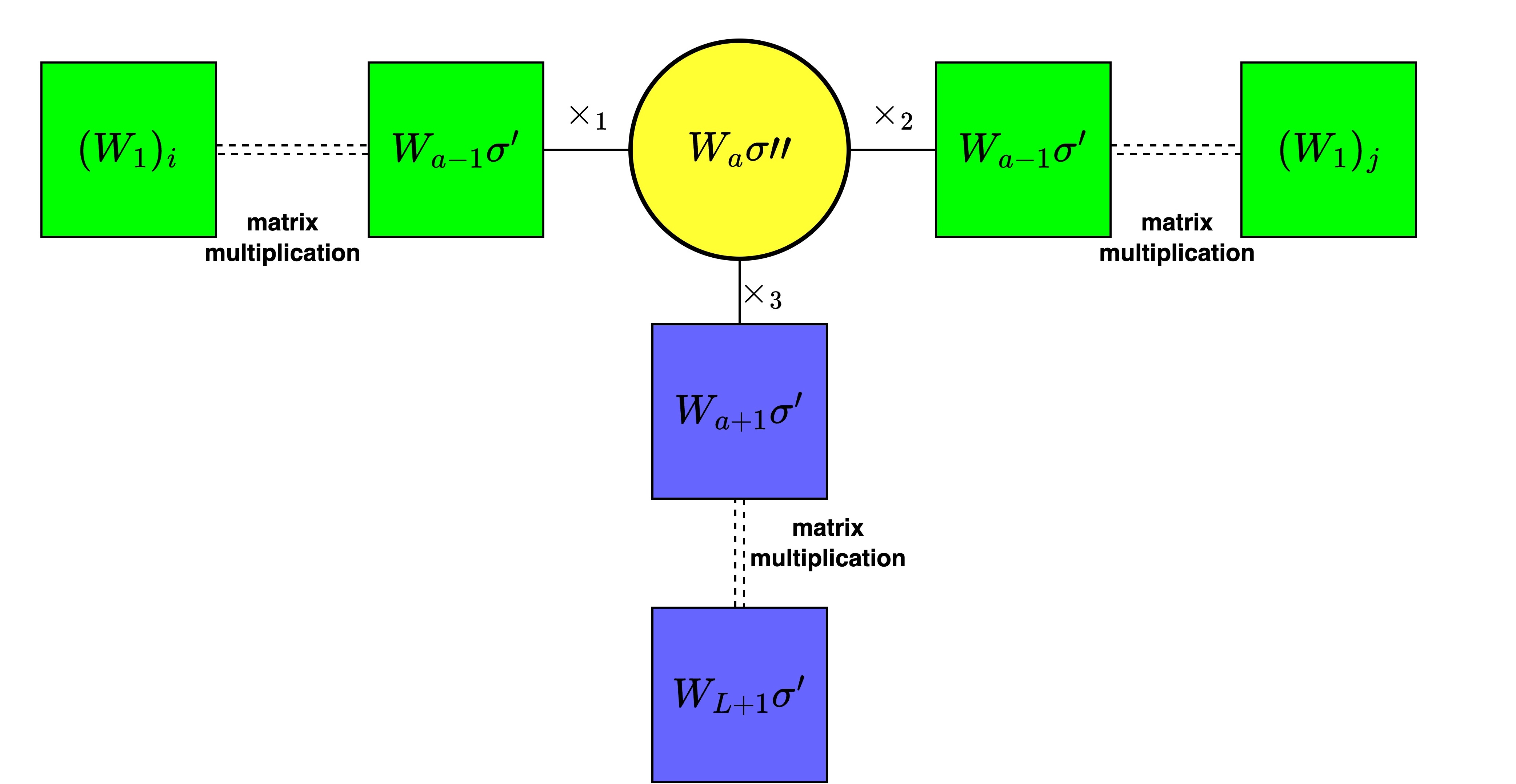}
				\caption{The structure of $\lambda_a(\vx)$.}
				\label{lambda}
			\end{figure}

Note that $\sigma$ can either be the ReLU or the ReLU squared function. If $\sigma(x)$ is the ReLU function, then its second derivative, $\sigma''(x)$, will be the delta distribution, which is defined as $\sigma''(x)=0$ for $x\neq0$ and $\sigma''(0)=\infty$, with $\int_\mathbb{R}\sigma''(x)\,\D x=1$. This delta function serves as the generalization function (a distribution. However, since we assume that $\phi\in W^2([a,b]^d)$, any delta functions must either disappear or be omitted in the $\lambda_a(\vx)$ equation.

Denote $W_i$ as the number of parameters in $\vW_i,\vb_i$, i.e., $W_i=N_iN_{i-1}+N_i$. Let $\vx\in\sR^d$ be an input and $\vtheta\in\sR^W$ be a parameter vector in $\psi$. We denote the output of $\psi$ with input $\vx$ and parameter vector $\vtheta$ as $f(\vx,\vtheta)$. For fixed $\vx_1,\vx_2,\ldots,\vx_m$ in $\sR^d$, we aim to bound\begin{align}
				K:=\left|\{\left(\sgn(f(\vx_1,\vtheta)),\ldots,\sgn(f(\vx_m,\vtheta))\right):\vtheta\in\sR^W\}\right|.
			\end{align}
			
The proof is follow our previous work \cite{yang2023nearly}. For any partition $\fS=\{P_1,P_2,\ldots,P_T\}$ of the parameter domain $\sR^W$, we have \[K\le \sum_{i=1}^T\left|\{\left(\sgn(f(\vx_1,\vtheta)),\ldots,\sgn(f(\vx_m,\vtheta))\right):\vtheta\in P_i\}\right|.\] We choose the partition such that within each region $P_i$, the functions $f(\vx_j,\cdot)$ are all fixed polynomials of bounded degree. This allows us to bound each term in the sum using Lemma \ref{bounded}.

  We define a sequence of sets of functions $\{\sF_j\}_{j=0}^L$ with respect to parameters $\vtheta\in\sR^W$:\begin{align}
				\sF_0:=\cup_{j=1}^m\{&(\vW_1)_1,(\vW_1)_2,\ldots,(\vW_1)_d,\vW_1\vx_j+\vb_1\}\notag\\
				\sF_1:=\cup_{j=1}^m\{&\vW_2\sigma''(\vW_1\vx_j+\vb_1),\vW_2\sigma'(\vW_1\vx_j+\vb_1),\vW_2\sigma(\vW_1\vx_j+\vb_1)+\vb_2\}\cup \sF_0\notag\\\sF_2:=\cup_{j=1}^m\{&\vW_3\sigma''(\vW_2\sigma(\vW_1\vx_j+\vb_1)+\vb_2),\vW_3\sigma'(\vW_2\sigma(\vW_1\vx_j+\vb_1)+\vb_2),\notag\\&\vW_3\sigma(\vW_2\sigma(\vW_1\vx_j+\vb_1)+\vb_2)+\vb_3\}\cup \sF_1\notag\\&\vdots\notag\\\sF_L:=\cup_{j=1}^m\{&\vW_{L+1}\sigma'(\vW_{L}\sigma_1(\ldots\sigma_1(\vW_1\vx_j+\vb_1)\ldots)+\vb_L),\notag\\&\vW_{L+1}\sigma''(\vW_{L}\sigma_1(\ldots\sigma_1(\vW_1\vx_j+\vb_1)\ldots)+\vb_L)\}\cup \sF_{L-1}.\notag
			\end{align}

   The partition of $\sR^W$ is constructed layer by layer through successive refinements denoted by $\fS_0,\fS_1,\ldots,\fS_L$. We denote $L^*=L-C\log_2 L$. These refinements possess the following properties:

            \textbf{1}. We have $|\fS_0|=1$, and for each $n=1,\ldots,L$, we have \[\frac{|\fS_n|}{|\fS_{n-1}|}\le 2\left(\frac{2em(1+(n-1)2^{\max\{0,n-1-L^*\}})N_n}{\sum_{i=1}^nW_i}\right)^{\sum_{i=1}^nW_i}.\]
			
			\textbf{2}. For each $n=0,\ldots,L^*$, each element $S$ of $\fS_{n}$, when $\vtheta$ varies in $S$, the output of each term in $\sF_n$ is a fixed polynomial function in $\sum_{i=1}^nW_i$ variables of $\vtheta$, with a total degree no more than $1+n2^{\max\{0,n-L^*\}}$.
			
			We define $\fS_0=\{\sR^W\}$, which satisfies properties 1,2 above, since $\vW_1\vx_j+\vb_1$ and $(\vW_1)_i$ for all $i=1,\ldots,d$ are affine functions of $\vW_1,\vb_1$.
			
			For each $n=0,\ldots,L$, to define $\fS_n$, we use the last term of $\sF_{n-1}$ as inputs for the new terms in $\sF_n$. All elements in $\sF_n$ except the $\sF_n\backslash \sF_{n-1}$ are fixed polynomial functions in $W_n$ variables of $\vtheta$, with a total degree no greater than $1+(n-1)2^{\max\{0,n-1-L^*\}}$ when $\vtheta$ varies in $S\in\fS_n$. This is because $\fS_n$ is a finer partition than $\fS_{n-1}$.
			
			We denote $p_{\vx_j,n-1,S,k}(\vtheta)$ as the output of the $k$-th node in the last term of $\sF_{n-1}$ in response to $\vx_j$ when $\vtheta\in S$. The collection of polynomials \[\{p_{\vx_j,n-1,S,k}(\vtheta): j=1,\ldots,m,~k=1,\ldots,N_n\}\]can attain at most $2\left(\frac{2em(1+(n-1)2^{\max\{0,n-1-L^*\}})N_n}{\sum_{i=1}^nW_i}\right)^{\sum_{i=1}^nW_i}$ distinct sign patterns when $\vtheta\in S$ due to Lemma \ref{bounded} for sufficiently large $m$. Therefore, we can divide $S$ into \[2\left(\frac{2em(1+(n-1)2^{\max\{0,n-1-L^*\}})N_n}{\sum_{i=1}^nW_i}\right)^{\sum_{i=1}^nW_i}\] parts, each having the property that $p_{\vx_j,n-1,S,k}(\vtheta)$ does not change sign within the subregion. By performing this for all $S\in\fS_{n-1}$, we obtain the desired partition $\fS_n$. This division ensures that the required property 1 is satisfied.
			
			Additionally, since the input to terms in $\sF_n\backslash\sF_{n-1}$ is $p_{\vx_j,n-1,S,k}(\vtheta)$, and we have shown that the sign of this input will not change in each region of $\fS_n$, it follows that the output of the terms in $\sF_n\backslash\sF_{n-1}$ is also a polynomial without breakpoints in each element of $\fS_n$, therefore, the required property 2 is satisfied.
Due to the structure of $\lambda$, it is a polynomial function in $\sum_{i=1}^{L+1}W_i$ variables of $\vtheta\in S\in\fS_L$, of total degree no more than \[d_2:=2\sum_{n=0}^L(1+n2^{\max\{0,n-L^*\}})=2L+2+L^C(L-1).\] Therefore, for each $S\in\fS_L$ we have \[\left|\{\left(\sgn(f(\vx_1,\vtheta)),\ldots,\sgn(f(\vx_m,\vtheta))\right):\vtheta\in S\}\right|\le 2\left(2emd_2/\sum_{i=1}^{L+1}W_i\right)^{\sum_{i=1}^{L+1}W_i}.\] Then 
			 \begin{align}
				K\le& 2\left(2emd_2/\sum_{i=1}^{L+1}W_i\right)^{\sum_{i=1}^{L+1}W_i}\cdot  \prod_{n=1}^{L}2\left(\frac{2em(1+(n-1)2^{\max\{0,n-1-L^*\}})N_n}{\sum_{i=1}^nW_i}\right)^{\sum_{i=1}^nW_i}\notag\\\le& 4\prod_{n=1}^{L+1}2\left(\frac{2em(1+(n-1)2^{\max\{0,n-1-L^*\}})N_n}{\sum_{i=1}^nW_i}\right)^{\sum_{i=1}^nW_i}\notag\\\le& 2^{L+3} \left(\frac{4em(L^{C+3})N}{U}\right)^{U},\notag
			\end{align}
   where $U:=\sum_{n=1}^{L+1}\sum_{i=1}^nW_i=\fO(N^2L^2)$, $N$ is the width of the network, and the last inequality is due to weighted AM-GM. For the definition of the VC-dimension, we have \begin{equation}
				2^{\operatorname{VCdim}(D^2\Phi)}\le 2^{L+3} \left(\frac{4e\operatorname{VCdim}(D^2\Phi)L^{C+3}N}{U}\right)^{U}.\notag
			\end{equation}
   
  Due to Lemma \ref{inequality}, we obtain that\begin{equation}
				\operatorname{VCdim}(D^2\Phi)=\fO(N^2L^2\log_2 L\log_2 N)\notag
			\end{equation}since $U=\fO(N^2L^2)$.
   \end{proof}
Our method for estimating the VC-dimension is also suitable for the higher-order derivatives of DNNs. The reason for this is that the partition of $\fS$ and the definition of the sets of function sets $\{\sF_j\}_{j=0}^L$ do not depend on the order of the derivatives in our method. Therefore, the VC-dimension of the higher-order derivatives of DSRNs is still $\fO(N^2L^2\log_2 L\log_2 N)$.

{\begin{remark}
   The estimation of the VC-dimension of higher-order derivatives of deep neural networks can follow our proof above. Based on this proof, we can notice that if the square of ReLU can appear at any position instead of just the last several layers like in DSRN, the upper bound of the VC-dimension will be \(\fO(N^2L^3\log N\log L)\) or \(\fO(N^3L^2\log N\log L)\) shown in \cite{bartlett2019nearly}, which is much more complex than DSRN. This will make the proof of optimality invalid. Furthermore, this can be one reason to show that deep neural networks with only the square of ReLU as the activation function are more difficult to train since a larger VC-dimension means increased complexity of the space.
\end{remark}}

			\begin{proof}[Proof of Theorem \ref{Optimality}]
				By combining Theorem \ref{vcdim} with the result presented in \cite[Theorem 5]{siegel2022optimal}, we can establish the proof for Theorem \ref{Optimality}.
			\end{proof}

   Furthermore, through the combination of Theorem \ref{main2} with the findings outlined in \cite[Theorem 5]{siegel2022optimal}, we can establish the following corollary. This corollary serves to demonstrate the optimality stated in Theorem \ref{vcdim}.
   \begin{corollary}\label{vcdim_opt}
			For any $d\in\sN_+$, $C,J_0,\varepsilon>0$, there exists $N,L\in\sN$ with $NL\ge J_0$ such that	\begin{equation}
				\operatorname{VCdim}(D^2\Phi)> C N^{2-\varepsilon}L^{2-\varepsilon},\end{equation}where $D^2\Phi$ is defined in Theorem \ref{vcdim}.
		\end{corollary}

 Note that in this section, we have proven the optimality of the approximation of DSRNs. However, if we allow the square of ReLU to appear in any layer of DNNs, the VC-dimension of derivatives of such DNNs will become $\fO(L^3N^2)$, as also shown in \cite{bartlett2019nearly}. In other words, we cannot prove the optimality of $\sigma_2$-NNs based on such a bound of VC-dimension. Nevertheless, we believe that the approximation rate $\fO(N^{-2(n-2)/d}L^{-2(n-2)/d})$ may not be improved to $\fO(N^{-2(n-2)/d}L^{-3(n-2)/d})$ for $\sigma_2$-NNs, as ReLU can efficiently approximate the square of ReLU \cite{yarotsky2017error}, and the optimality of ReLU DNNs has been proven in \cite{lu2021deep,yang2023nearly}. Obtaining optimal bounds for $\sigma_2$-NNs will be a future research direction.

\section{Generalization Analysis in Sobolev Spaces}
Let $\rho$ be a Borel probability measure on
\(
  \mathcal Z := \Omega \times \mathcal Y
\)
with $\mathcal Y\subset\mathbb R$.  
Assume that the $\vx$–marginal of $\rho$ is the uniform distribution on $\Omega$, and that the conditional law $\rho(\mathrm d y \,|\,\vx)$ satisfies  
\[
  \fL f_{\rho}(\vx)
  \;=\;
  \int_{\mathcal Y} y \,\rho(\mathrm d y \,|\,\vx),
  \qquad
  f_{\rho} \in W^{n,\infty}(\Omega),
\]
where $\fL$ satisfies Assumption~\ref{assump:L}.

\begin{assumption}\label{assump:L}
$\fL$ is a second-order linear differential operator, and there exists a constant $C>0$ such that, for any $f_1,f_2 \in H^2(\Omega)$,
\[
  \|\fL f_1 - \fL f_2\|_{L^2(\Omega)} \le C \, \|f_1 - f_2\|_{H^2(\Omega)}.
\]
\end{assumption}

Draw
\(
  \mathcal S
  =\bigl\{(\vx_{j},y_{j})\bigr\}_{j=1}^{M}\subset\mathcal Z^{M}
\)
i.i.d.\ according to $\rho$.  
Define the energy loss
\[
  \mathcal E(f)
  =\int_{\mathcal Z}\bigl(\fL f(\vx)-y\bigr)^{2}\,\rho(\mathrm d\vx,\mathrm d y),
  \qquad
  \mathcal E_{\mathcal S}(f)
  =\frac1M\sum_{j=1}^{M}\bigl(\fL f(\vx_{j})-y_{j}\bigr)^{2}.
\]
Motivated by Theorem~\ref{main2}, set 

\begin{equation}
  \mathcal F_{N,L,D_1,D_2,D_3,B}
  :=\Bigl\{
      \phi
      \,\bigm|\,
      \phi\in\fN_{D_1,D_2 L\log L}\text{ with width}\le D_3 N\log N, \|\phi\|_{W^{2,\infty}(\Omega)}\le B
    \Bigr\},\label{set}
\end{equation}     
and $\fL \mathcal F_{N,L,D_1,D_2,D_3,B}:=\{\fL f\mid f\in \mathcal F_{N,L,D_1,D_2,D_3,B}\}$.
Let
\[
  f_{\mathcal S,\mathcal F_{N,L,D_1,D_2,D_3,B}}
  :=\arg\min_{f\in\mathcal F_{N,L,D_1,D_2,D_3,B}}\mathcal E_{\mathcal S}(f)
\]
be the empirical risk minimizer over $\mathcal F_{N,L,D_1,D_2,D_3,B}$.  
In what follows we bound the generalization error
\[
  \mathbb E\bigl\|\fL f_{\mathcal S,\mathcal F_{N,L,D_1,D_2,D_3,B}}-\fL f_{\rho}\bigr\|_{L^2(\Omega)},
\]
where the expectation is taken with respect to the sampling of $\mathcal S$. 

The principal theorem of this section is stated below:\begin{theorem}\label{general thm}
Let \(f_{\rho} \in W^{n,\infty}(\Omega)\) and suppose Assumption~\ref{assump:L} holds.
Choose integers \(N, L, d \in \mathbb{N}_{+}\) satisfying \(\log N \ge d(\log n + \log d)\).
Assume there exist constants \(B, D_{1}, D_{2}, D_{3} > 0\) such that 
\(\|\fL f_{\rho}\|_{W^{2,\infty}(\Omega)} \le B\), and 
\[
\inf_{f \in \mathcal{F}_{N,L,D_{1},D_{2},D_{3},B}}
  \|f - f_{\rho}\|_{W^{2,\infty}(\Omega)}
\]
achieves the approximation rate stated in Theorem~\ref{main2}.
Moreover, suppose the responses satisfy \(\mathcal{Y} \subset [-L, L]\) almost surely.  
Then there exists \(M_{0} > 0\) such that for all \(M \ge M_{0}\),
\[
\mathbb{E}\bigl\|
    \fL f_{\mathcal{S},\mathcal{F}_{N,L,D_{1},D_{2},D_{3},B}}
    - \fL f_{\rho}
\bigr\|_{L^{2}(\Omega)}
\;\le\;
C\,M^{-\frac{2(n-2)}{d+2(n-2)}}(\log M)^{\frac{6(n-2)}{d+2(n-2)}},
\]
up to lower-order \(\log\log M\) factors, where \(C>0\) depends only polynomially on 
\(d\), \(B\), and the constants \(D_{1}, D_{2}, D_{3}\).
\end{theorem}

Compared with the bound in \cite{suzukiadaptivity}, our estimate sharpens the logarithmic factor while preserving the same polynomial dependence.  
In the notation of \cite[Theorem 2]{suzukiadaptivity}, this corresponds to choosing \(s=n-2\), which is natural in our setting because the function has Sobolev regularity \(n\) and the measure norm is \(2\); the difference between these two values is therefore exactly \(n-2\).  
Consequently, the rate we obtain is nearly optimal. Our generalization analysis departs from the frameworks of \cite{suzukiadaptivity,schmidt2020nonparametric}, which assume uniformly bounded network parameters.  
In our construction, certain parameters grow exponentially in~$d$ for two principal reasons based on Propositions \ref{step} and \ref{point}
These effects necessitate a different analytical approach to control the estimation error.

Instead, we establish the generalization bound through the \emph{pseudo-dimension} (Definition \ref{Pse}) of the network’s  second order derivatives.  
To connect the generalization error with the pseudo-dimension we invoke Lemma~\ref{connect} below, which rests on covering-number estimates.  
We therefore start by recalling the relevant notion of covering numbers.
\begin{definition}[covering number \cite{anthony1999neural}]
				Let $(V,\|\cdot\|)$ be a normed space, and $\Theta\subset V$. $\{V_1,V_2,\ldots,V_n\}$ is an $\varepsilon$-covering of $\Theta$ if $\Theta\subset \cup_{i=1}^nB_{\varepsilon,\|\cdot\|}(V_i)$. The covering number $\fN(\varepsilon,\Theta,\|\cdot\|)$ is defined as \[\fN(\varepsilon,\Theta,\|\cdot\|):=\min \{n: \exists \varepsilon \text {-covering over } \Theta \text { of size } n\} \text {. }\]
			\end{definition}
			
			\begin{definition}[Uniform covering number \cite{anthony1999neural}]\label{uniform}
				{Suppose the $\fF$ is a class of functions from $\fX$ to $\sR$.} Given $n$ samples $\vZ_n=(z_1,\ldots,z_n)\in\fX^n$, define \[\fF|_{\vZ_n}=\{(u(z_1),\ldots,u(z_n)):u\in\fF\}.\]The uniform covering number $\fN(\varepsilon,\fF,n)$ is defined as \[\fN(\varepsilon,\fF,n)=\max_{\vZ_n\in\fX^n}\fN\left(\varepsilon, \fF|_{\vZ_n},\|\cdot\|_{\infty}\right),\]where $\fN\left(\varepsilon, \fF|_{\vZ_n},\|\cdot\|_{\infty}\right)$ denotes the $\varepsilon$-covering number of $\fF|_{\vZ_n}$ w.r.t the $l^\infty$-norm on $\vZ_n$ defined as $\|f\|_{\infty}=\sup_{\vz_i\in\vZ_n}|f(z_i)|$.
			\end{definition}

\begin{lemma}[\cite{gyorfi2002distribution}, Theorem 11.4]\label{connect}
				Let $M \in \sN$, and assume that $\|g\|_{L^{\infty}(\Omega)}\le B$ and $\fY\in[-B,B]$ are almost surely for some $B \geq 1$. Let $\fF$ be a set of functions from $\Omega$ to $[-B, B]$. Then for any $0<\delta \leq 1 / 2$ and $\alpha, \beta>0$,
\[
\begin{aligned}
& \mathbb{P}\left\{\exists \phi \in \fF:\left\|f-\phi\right\|_{L^2}^2-\left(\mathcal E_{\mathcal S}(\phi)-\mathcal E_{\mathcal S}\left(g\right)\right) \geq \epsilon\left(\alpha+\beta+\left\|g-\phi\right\|_{L^2}^2\right)\right\} \\
& \leq 14 \fN\left(\frac{\beta \delta}{20 B}, \fF, M\right) \exp \left(-\frac{\delta^2(1-\delta) \alpha M}{214(1+\delta) B^4}\right).
\end{aligned}
\]
			\end{lemma}
In the ordinary Lemma of \cite{gyorfi2002distribution}, the covering number of \(l^1\) is used to bound the inequality, which is stronger than our lemma here. Since \(l^\infty\) is sufficient to obtain our results and to use Lemma \ref{cover dim} directly, we utilize the covering number under the \(l^\infty\)-norm here.

With the preparations, we are ready to bound the generalization error by combining the approximation error with a covering-number estimate, as stated in the following proposition:
\begin{proposition}\label{prop:connect}
    Let $M \in \sN$, and assume that $\|\fL f_\rho\|_{L^{\infty}(\Omega)}\le B$ and $\fY\in[-L,L]$ are almost surely for some $B \geq 1$. Then we have that \begin{align}
        &\mathbb{E}\|\fL f_{\fS,\fF_{N,L,D_1,D_2,D_3,B}}-\fL f_\rho\|_{L^2(\Omega)}^2\le 2 \inf _{f \in \fF_{N,L,D_1,D_2,D_3,B}}\left\|\fL f-\fL f_\rho\right\|_{L^2(\Omega)}^2\notag\\+&\frac{5136 B^4}{M}\left\{\log \left(14 \ \mathcal{N}\left(\frac{1}{80 B M}, \fL\fF_{N,L,D_1,D_2,D_3,B},M\right)\right)+1\right\}.\notag
    \end{align}
\end{proposition}

\begin{proof}
    To simplicity notations, we denote $\fF_{N,L,D_1,D_2,D_3,B}$ as $\fF_*$ in this proof. First, we divide \(\|\fL f_{\fS,\fF_*}-\fL f_\rho\|_{L^2(\Omega)}^2\) into two parts \begin{align*}&\|\fL f_{\fS,\fF_*}-\fL f_\rho\|_{L^2(\Omega)}^2=2(\fE_{\fS}(f_{\fS,\fF_*})-\fE_{\fS}(f_\rho))\notag\\+&\underbrace{\|\fL f_{\fS,\fF_*}-\fL f_\rho\|_{L^2(\Omega)}^2-2(\fE_{\fS}(f_{\fS,\fF_*})-\fE_{\fS}(f_\rho))}_{\mathbf{A}_0}.\end{align*}For the first part, we have that \[\begin{aligned}
        &2\mathbb{E}(\fE_{\fS}(f_{\fS,\fF_*})-\fE_{\fS}(f_\rho))\le2\mathbb{E}\inf_{f\in\fF_*}(\fE_{\fS}(f)-\fE_{\fS}(f_\rho))\le 2\inf_{f\in\fF_*}\mathbb{E}(\fE_{\fS}(f)-\fE_{\fS}(f_\rho))\\=&2\inf_{f\in\fF_*}\mathbb{E}\left(\frac{1}{M}\sum_{j=1}^M|\fL f(\vx_i)-\fL f_\rho(\vx_i)|^2\right)\le 2 \inf _{f \in \fF_*}\left\|\fL f-\fL f_\rho\right\|_{L^2(\Omega)}^2,\notag
    \end{aligned}\]where the first equality is use to 
      \[\mathbb{E}_{\vx_i}(\fL f(\vx_i)-\fL f_\rho(\vx_i))\mathbb{E}_{y_i}[(\fL f_\rho(\vx_i)-y_i)|\vx_i]=\mathbb{E}_{\vx_i}(\fL f(\vx_i)-\fL f_\rho(\vx_i))\cdot 0=0,\]and last inequality is due to Assumption \ref{assump:L}. As for $\mathbf{A}_0$, we have \begin{align}
    &\quad\mathbb{P}(\mathbf{A}_0\ge \epsilon)\notag\\&=\mathbb{P}\left( 2\|\fL f_{\fS,\fF_*}-\fL f_\rho\|_{L^2(\Omega)}^2-2(\fE_{\fS}(f_{\fS,\fF_*})-\fE_{\fS}(f_\rho))\ge\epsilon+\|\fL f_{\fS,\fF_*}-\fL f_\rho\|_{L^2(\Omega)}^2\right)\notag\\&=\mathbb{P}\left( \|\fL f_{\fS,\fF_*}-\fL f_\rho\|_{L^2(\Omega)}^2-(\fE_{\fS}(f_{\fS,\fF_*})-\fE_{\fS}(f_\rho))\ge\frac{1}{2}\left(\epsilon+\|\fL f_{\fS,\fF_*}-\fL f_\rho\|_{L^2(\Omega)}^2\right)\right)\notag\\&\le 14 \fN\left(\frac{ \epsilon}{80 B}, \fL\fF_*,M\right) \exp \left(-\frac{\epsilon M}{5136 B^4}\right),\notag
\end{align}where the last inequality is due to Lemma \ref{connect} for choosing $\alpha=\beta=\frac{1}{2}\epsilon$ and $\delta=\frac{1}{2}$.

Therefore, we have that \[\begin{aligned}
    \mathbb{E} \mathbf{A}_0&\le \int_{0}^\infty \mathbb{P}(\mathbf{A}_0\ge t)\,\D t\le\epsilon+\int_{\epsilon}^\infty \mathbb{P}(\mathbf{A}_0\ge t)\,\D y\\&\le\epsilon+\int_{\epsilon}^\infty 14 \fN\left(\frac{ \epsilon}{80 B}, \fL\fF_*,M\right) \exp \left(-\frac{t M}{5136 B^4}\right)\,\D t\notag
\end{aligned}\]

By the direct calculation, we have \begin{align}
    &\int_{\epsilon}^\infty 14 \fN\left(\frac{ \epsilon}{80 B}, \fL\fF_*,M\right) \exp \left(-\frac{t M}{5136 B^4}\right)\,\D t\notag\\\le& 14 \fN\left(\frac{ \epsilon}{80 B}, \fL\fF_*,M\right) \frac{5136 B^4}{M}\exp \left(-\frac{\epsilon M}{5136 B^4}\right).\notag
\end{align} Set \[\epsilon=\frac{5136 B^4}{M}\log\left(14 \fN\left(\frac{ 1}{80 BM}, \fL\fF_*,M\right) \right)\ge \frac{1}{M}\] and we have \[\mathbb{E} \mathbf{A}_0\le \frac{5136 B^4}{M}\left[\log\left(14 \fN\left(\frac{ 1}{80 BM}, \fL\fF_*,M\right) \right)+1\right]\]

Hence we have \begin{align}
        &\mathbb{E}\|\fL f_{\fS,\fF_*}-\fL f_\rho\|_{L^2(\Omega)}^2\notag\\\le& \frac{5136 B^4}{M}\left[\log\left(14 \fN\left(\frac{ 1}{80 BM}, \fL\fF_*,M\right) \right)+1\right]+2\inf_{f\in\fF_*}\left\|\fL f-\fL f_\rho\right\|_{L^2(\Omega)}^2.\notag
    \end{align}
\end{proof}

We bound the covering numbers via the pseudo-dimension, defined below.
\begin{definition}
   [pseudo-dimension \cite{pollard1990empirical}]\label{Pse}
		Let $\fF$ be a class of functions from $\fX$ to $\sR$. The pseudo-dimension of $\fF$, denoted by $\operatorname{Pdim}(\fF)$, is the largest integer $m$ for which there exists $(x_1,x_2,\ldots,x_m,y_1,y_2,\ldots,y_m)\in\fX^m\times \sR^m$ such that for any $(b_1,\ldots,b_m)\in\{0,1\}^m$ there is $f\in\fF$ such that $\forall i: f\left(x_i\right)>y_i \Longleftrightarrow b_i=1.$    
   \end{definition}

\begin{lemma}[\cite{anthony1999neural}, Theorem 12.2]\label{cover dim}
				Let $\fF$ be a class of functions from $\fX$ to $[-B,B]$. For any $\varepsilon>0$, we have \[\fN(\varepsilon,\fF,n)\le \left(\frac{2\mathrm{e}nB}{\varepsilon\operatorname{Pdim}(\fF)}\right)^{\operatorname{Pdim}(\fF)}\] for $n\ge \operatorname{Pdim}(\fF)$.
			\end{lemma}

\begin{proposition}\label{prop:pdim}
		For any $N,L,d,D_i,B\in\sN_+$, there exists a constant $\bar{C}$ independent of $N,L$ such that	\begin{equation}
		\operatorname{Pdim}(\fL \mathcal F_{N,L,D_1,D_2,D_3,B})\le \bar{C} N^2L^2\log_2 L\log_2 N,\notag
  \end{equation}where $\fL \mathcal F_{N,L,D_1,D_2,D_3,B}$ is defined in \eqref{set} and $\fL$ satisfies Assumption \ref{assump:L}.
	\end{proposition}
			
	\begin{proof}
Denote
\[
\Phi_{\fN}
:=
\{\eta(\vx,y) : \eta(\vx,y) = \psi(\vx) - y,\;
\psi \in \fL \mathcal{F}_{N,L,D_1,D_2,D_3,B},\;
(\vx,y) \in \mathbb{R}^{d+1} \}.
\]
By the definitions of VC-dimension and pseudo-dimension, we have
\begin{equation}
\operatorname{Pdim}(\fL \mathcal{F}_{N,L,D_1,D_2,D_3,B})
\;\le\;
\operatorname{VCdim}(\Phi_{\fN}).
\end{equation}
For \(\operatorname{VCdim}(\Phi_{\fN})\), it can be bounded by 
\(\mathcal{O}(N^2 L^2 \log_2 L \log_2 N)\).
The proof is analogous to that for the estimate of
\(\operatorname{VCdim}(D^2\Phi)\) given in Theorem~\ref{vcdim}.
The form of \(\fL\) does not affect the proof since \(\fL\) is a linear operator;
for any fixed input \(\vx\), the output remains a piecewise polynomial with respect
to the parameter space. Therefore, the argument used for
\(\operatorname{VCdim}(D^2\Phi)\) in Theorem~\ref{vcdim} still applies here.
\end{proof}

\begin{proof}[Proof of Theorem \ref{general thm}]
    % Denote that \[\fF_*=\{s\mid s\text{ is a squared–ReLU network on $\Omega$ of width }3d^{2}
    %   \text{ and depth }\lfloor\log_{2}d\rfloor+2\}.\]Based on $\nabla \fF_{N,L,D_1,D_2,D_3,B}\subset \sum_{s=1}^m\nabla \fF_*$ and \cite[p.~21]{sen2018gentle}, we have that \[\fN\left(\frac{ 1}{80 BM}, \fL\fF_{N,L,D_1,D_2,D_3,B},M\right)\le\fN\left(\frac{ 1}{80 BMm}, \fL\fF_*,M\right)^m.\]Therefore, b
    Based on Lemma~\ref{cover dim}, Propositions~\ref{prop:connect} and~\ref{prop:pdim}, set $M_0=\operatorname{Pdim}(\fL\fF_{N,L,D_1,D_2,D_3,B})$, we have that for any $M\ge M_0$, \begin{align}
        &\mathbb{E}\|\fL f_{\fS,\fF_{N,L,D_1,D_2,D_3,B}}-\fL f_\rho\|_{L^2(\Omega)}^2\notag\\\le& \frac{5136 B^4}{M}\left[\log\left(14 \fN\left(\frac{ 1}{80 BM}, \fL\fF_{N,L,D_1,D_2,D_3,B},M\right) \right)+1\right]+C_*(NL)^{-\frac{4(n-2)}{d}}\notag\\\le &\frac{5136 B^4}{M}\left[\log\left( \fN\left(\frac{ 1}{80 BM}, \fL\fF_{N,L,D_1,D_2,D_3,B},M\right) \right)+4\right]+C_*(NL)^{-\frac{4(n-2)}{d}}\notag\\\le& \frac{5136 B^4}{M}\left[\operatorname{Pdim}(\fL\fF_{N,L,D_1,D_2,D_3,B})\log\left(160\mathrm{e}M^2B^2 \right)+4\right]+C_*(NL)^{-\frac{4(n-2)}{d}}\notag\\\le & \frac{5136 B^4}{M}\left[\bar{C}N^2L^2\log_2N\log_2L\log\left(160\mathrm{e}M^2B^2 \right)+4\right]+C_*(NL)^{-\frac{4(n-2)}{d}}\notag\\=& \fO\left(\frac{N^2L^2\log N\log L\log M}{M}+(NL)^{-\frac{4(n-2)}{d}}\right).\notag
    \end{align}We consider
\[
E_1 \;=\; \frac{N^2L^2\,(\log N)(\log L)\,\\log M}{M},
\qquad
E_2 \;=\; (NL)^{-\frac{4(n-2)}{d}}.
\]
Let \(P:=NL\) and assume \(N=\fO(L)\), so that \((\log N)(\log L)=\fO\left((\log P)^2\right)\).
Balancing the two terms, \(E_1=E_2\), yields
\[
P^{\,2+\alpha}(\log P)^2 \;=\; \fO\left(\frac{M}{\log M}\right),
\qquad
\alpha:=\frac{4(n-2)}{d}.
\]
Solving for the optimal \(P\) gives
\[
P^\star =\fO\left(\; M^{\frac{1}{2+\alpha}}\,(\log M)^{-\frac{3}{2+\alpha}}\right)
\quad\text{(ignoring lower-order \(\log\log M\) factors)}.
\]
At the balance point, the minimal error satisfies
\[
M^{-\frac{2(n-2)}{\,d+2(n-2)\,}}\,
(\log M)^{\frac{6(n-2)}{\,d+2(n-2)\,}}.
\]

\end{proof}
\section{Numerical Experiment}\label{num}
In this section, we provide some examples to validate that DSRN can approximate smooth functions in Sobolev spaces.
\subsection{General function approximation}
\begin{figure}[!th]
\centering
        \begin{subfigure}[t]{0.32\textwidth}
		\centering
		\includegraphics[width=0.99\textwidth]{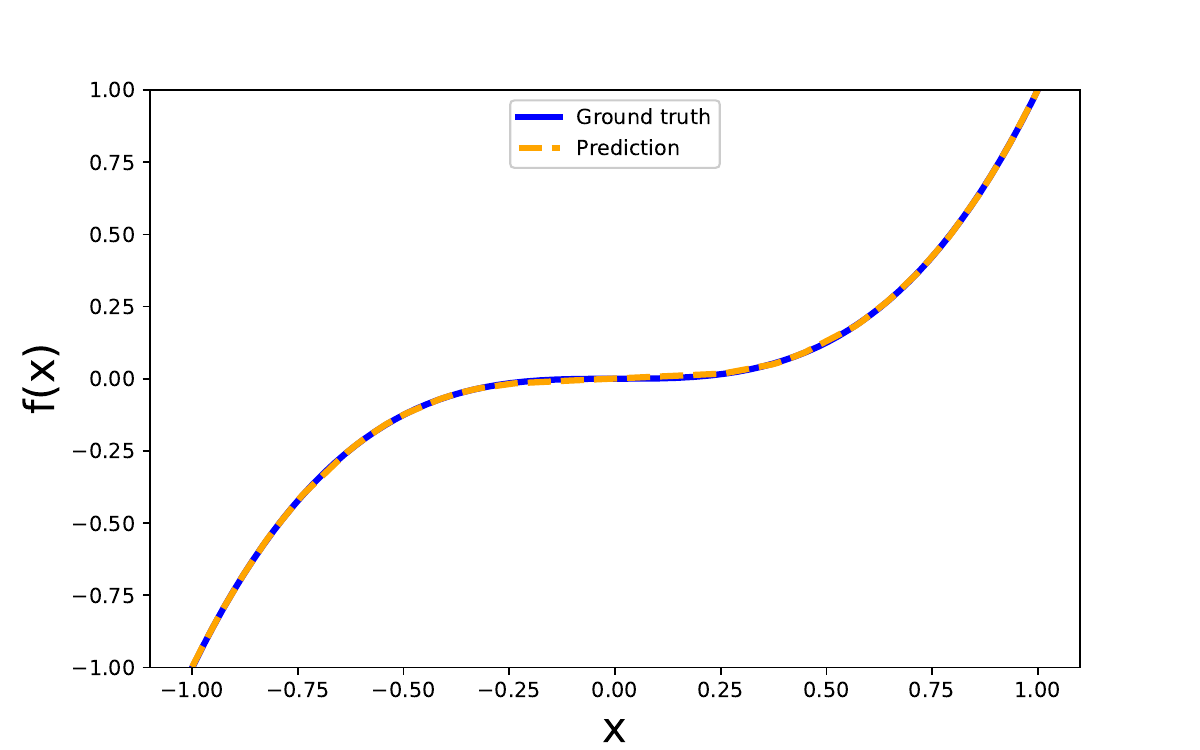}%V5
		\caption{ReLU: $f(x)=x^3$}
		\label{fig:relu-v}
	\end{subfigure}	
	\begin{subfigure}[t]{0.32\textwidth}
		\centering
		\includegraphics[width=0.99\textwidth]{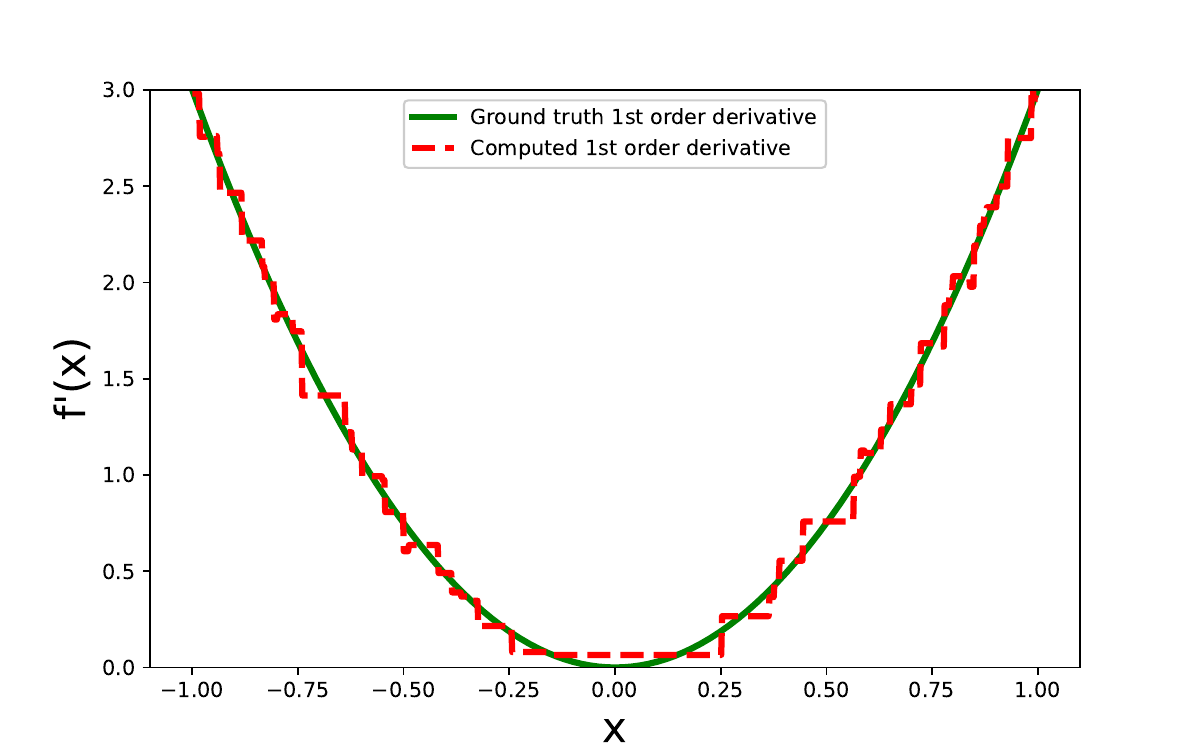}%V5
		\caption{ReLU: $f'(x)=3 x^2$}
		\label{fig:relu-g1}
	\end{subfigure}	
    \begin{subfigure}[t]{0.32\textwidth}
		\centering
		\includegraphics[width=0.99\textwidth]{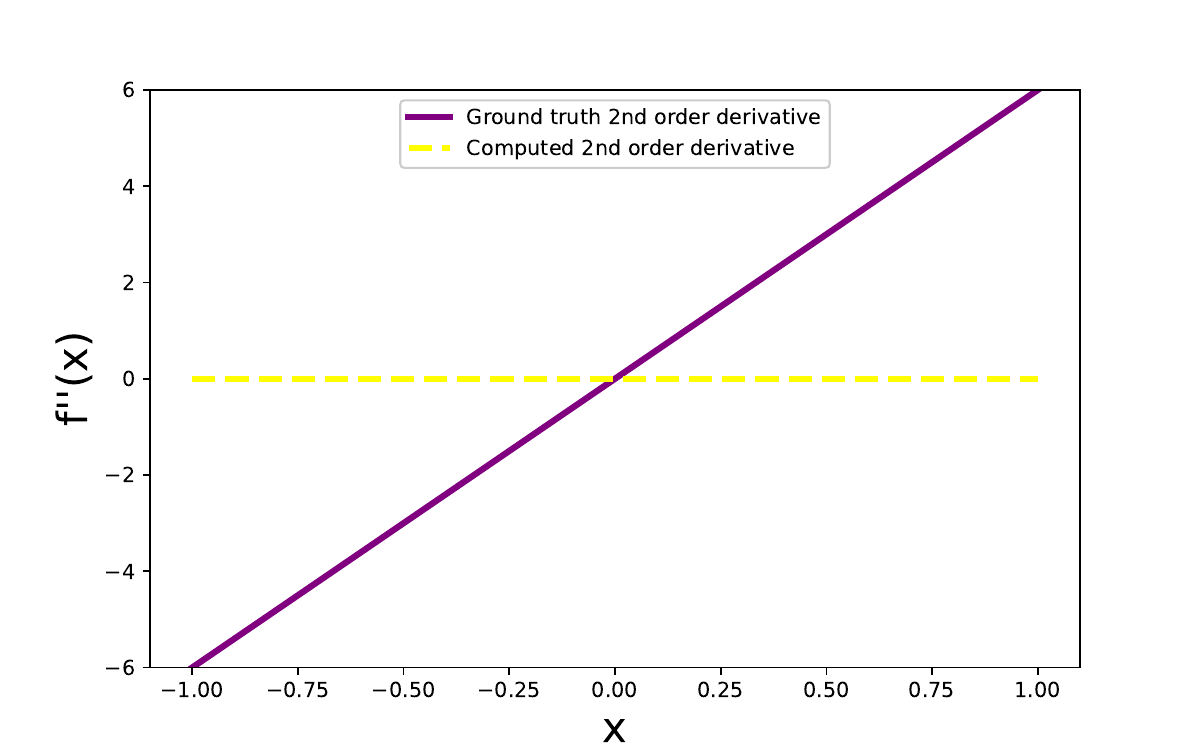}%V5
		\caption{ReLU: $f''(x)=6x$}
		\label{fig:relu-g2}
	\end{subfigure}	

	\begin{subfigure}[t]{0.32\textwidth}
		\centering
		\includegraphics[width=0.99\textwidth]{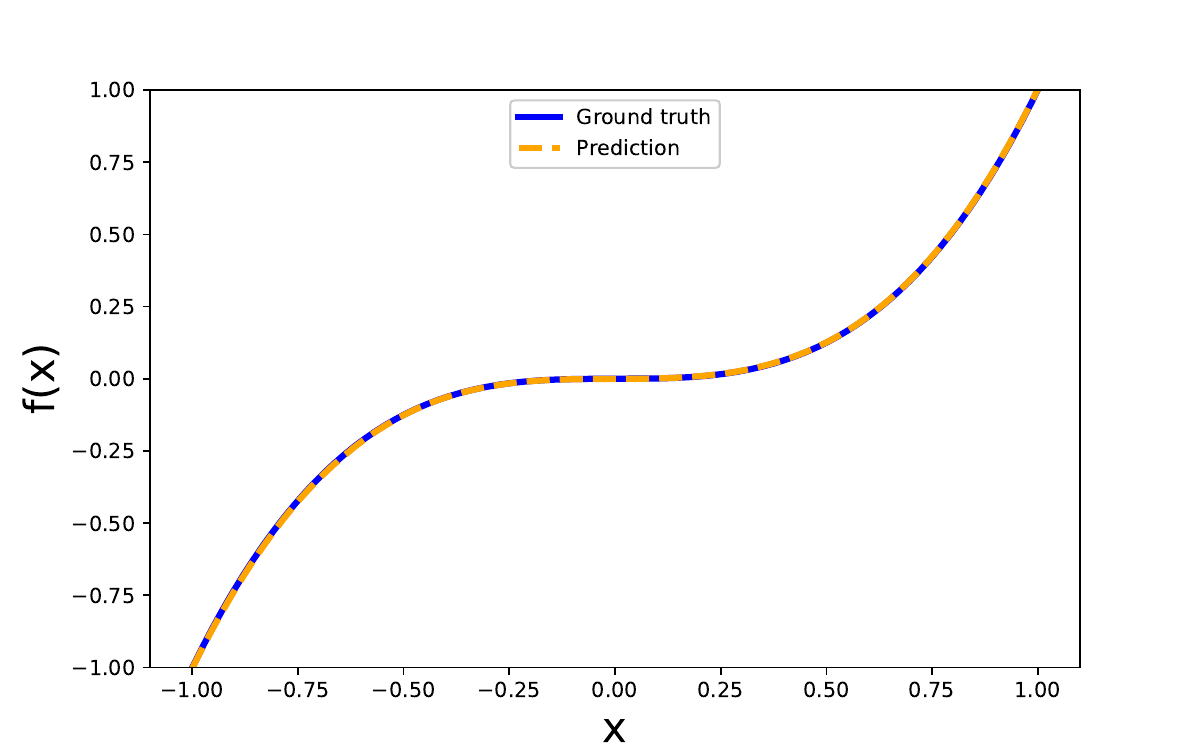}%V5
		\caption{DSRN: $f(x)=x^3$}
		\label{fig:dsrn-v}
	\end{subfigure}	
       \begin{subfigure}[t]{0.32\textwidth}
		\centering
		\includegraphics[width=0.99\textwidth]{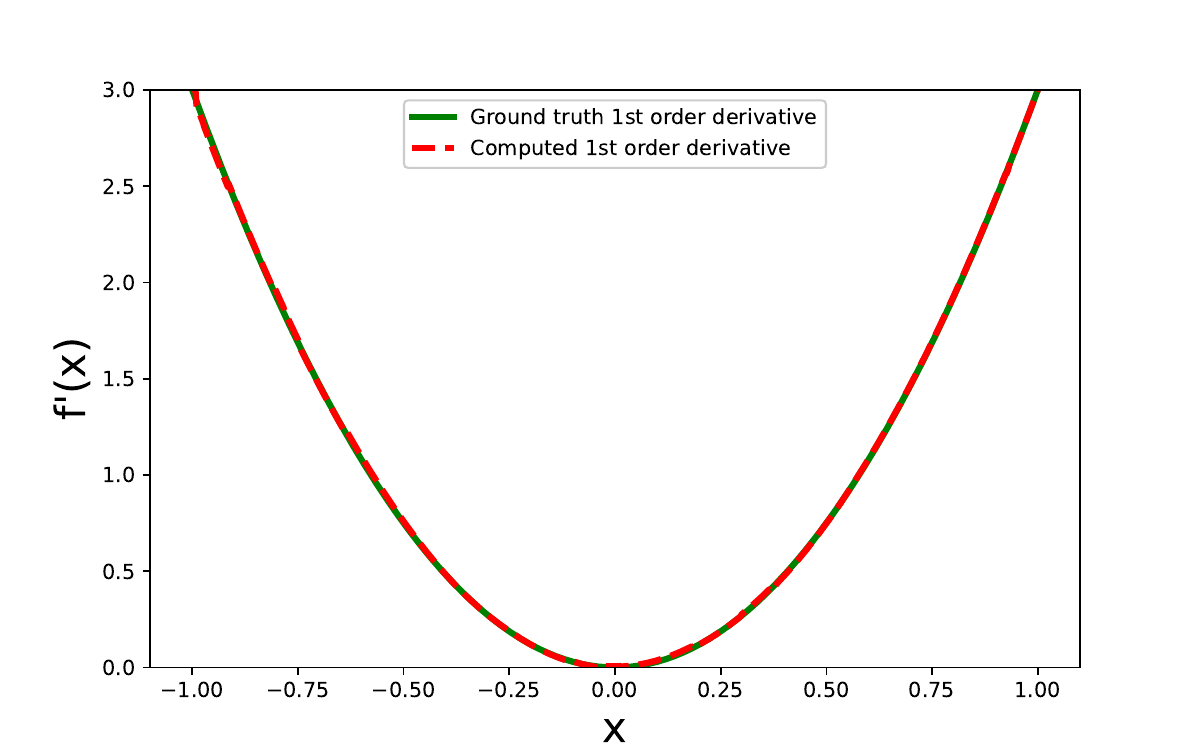}%V5
		\caption{DSRN: $f'(x)=3 x^2$}
		\label{fig:dsrn-g1}
	\end{subfigure}	
    \begin{subfigure}[t]{0.32\textwidth}
		\centering
		\includegraphics[width=0.99\textwidth]{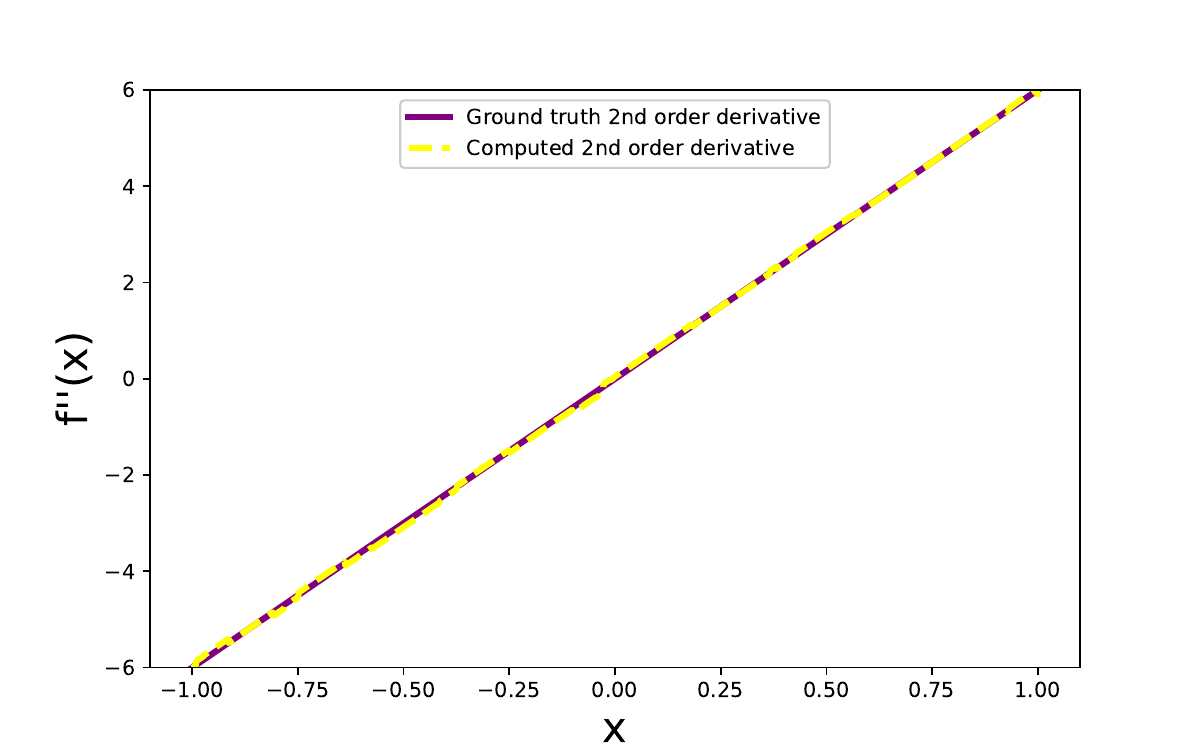}%V5
		\caption{DSRN: $f''(x)=6x$}
		\label{fig:dsrn-g2}
	\end{subfigure}	
	\caption{Comparison of function prediction and the first-order, second-order derivatives between true value, ReLU network and DSRN.}	
	\label{fig:xcubic}
\end{figure}

To highlight the advantage of Squared ReLU activation in DSRN for approximation in Sobolev norm, we compare it against a standard ReLU neural network in learning the cubic function $f(x)=x^3$ as well as its first- and second-order derivatives.
We separately train a 10-layer fully-connected  ReLU neural network, and a DSRN, where the last three layers use squared ReLU activations, to approximate the function $f(x)$. 
Fig.~\ref{fig:xcubic} presents the functions learned by the ReLU network (Fig.~\ref{fig:relu-v}) and DSRN (Fig.~\ref{fig:dsrn-v}),
the first-order derivative $f'(x)=3x^2$
estimated by the ReLU neural network (Fig.~\ref{fig:relu-g1}) and DSRN (Fig.~\ref{fig:dsrn-g1}),
and
the second-order derivative $f''(x)=6x$
estimated by the ReLU neural network (Fig.~\ref{fig:relu-g2}) and DSRN (Fig.~\ref{fig:dsrn-g2}).
The comparison results demonstrate that the DSRN successfully approximates both the cubic function and its derivatives, capturing smoothness in higher-order Sobolev norms.
In contrast, the ReLU neural network can only approximate the function up to the first-order derivative, but entirely fails to approximate the second-order derivative due to the vanishing curvature of ReLU.
Consequently, ReLU incurs significant error in higher-order Sobolev approximation, whereas DSRN achieves accurate approximation in higher-order Sobolev norm.

\FloatBarrier
\subsection{One-dimensional Poisson equation}
We solve the Poisson equation
\begin{equation}
    \Delta u(x) = -\pi^2 \sin(\pi x),  \quad  x \in \left[-1, 1\right],
\end{equation}
with boundary conditions
$u(-1)=u(1)=0$.
The exact solution is $u(x) = \sin(\pi x)$.
To be specific, we approximate the solution $u(x)$ by the neural network $\mathcal{NN}_{\vtheta}(x)$ and train the neural network parameters $\vtheta$ to minimize the mean squared error(MSE)
\[ \frac{1}{N} \sum_{i=1}^N 
\left[ 
\left(\frac{d^2}{d x^2}  \mathcal{NN}_{\vtheta} (x_i) + \pi^2 \sin(\pi x_i) \right)^2  
+
\left(\mathcal{NN}_{\vtheta} (x_i) -  y_i \right)^2  
\right].\]
The first term comes from PINN method and the second term is used for supervised leaning with $N$ training data $\left\{ (x_i,y_i),i=1,\dots,N | y_i=\sin(\pi x_i)\right\}$. 
We train 14-layer fully-connected networks to learn the solution $u(x)$.
For DSRN, the activation function of the last four hidden layers is the square of ReLU, and all other activation functions are ReLU.
For comparison, we also train neural networks with ReLU, and the square of ReLU as the activation function, respectively. 
We use the Adam optimizer with learning step $0.001$, weight decay $0.00001$ in the training.

\begin{figure}[!h]
\centering
        \begin{subfigure}[t]{0.45\textwidth}
		\centering
		\includegraphics[width=0.99\textwidth]{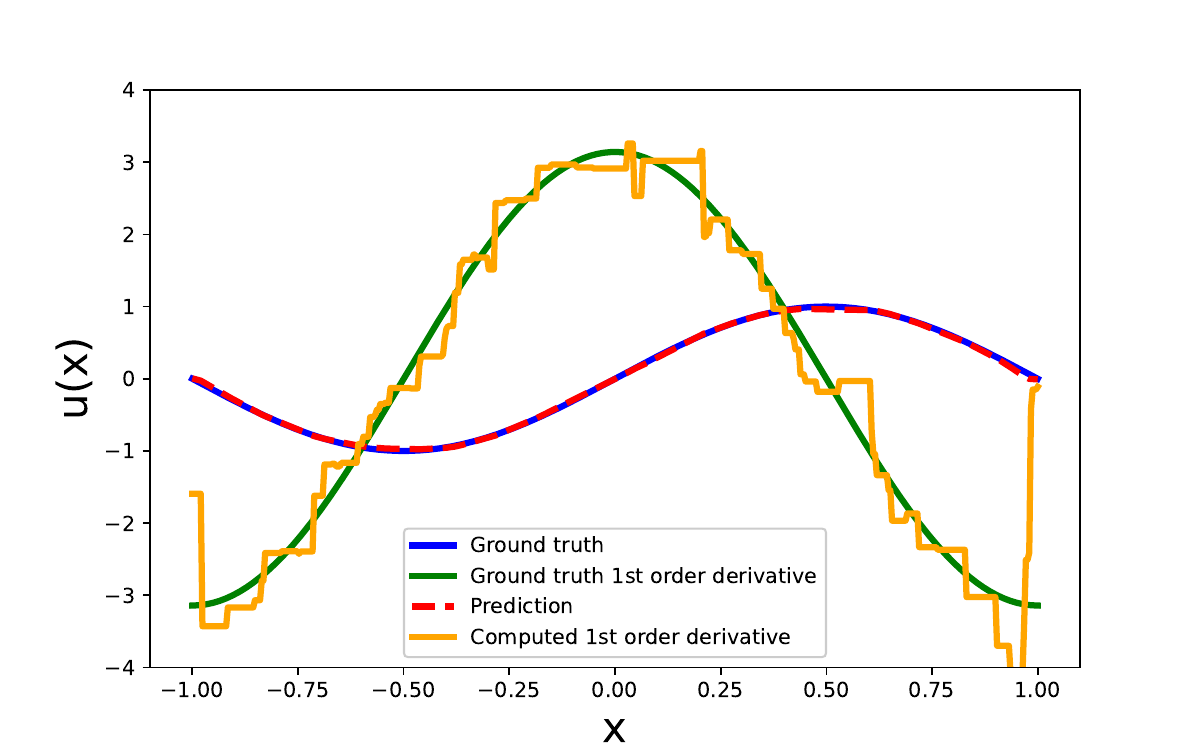}%V5
		\caption{ReLU}
		\label{fig:compare-a}
	\end{subfigure}	
	\begin{subfigure}[t]{0.45\textwidth}
		\centering
		\includegraphics[width=0.99\textwidth]{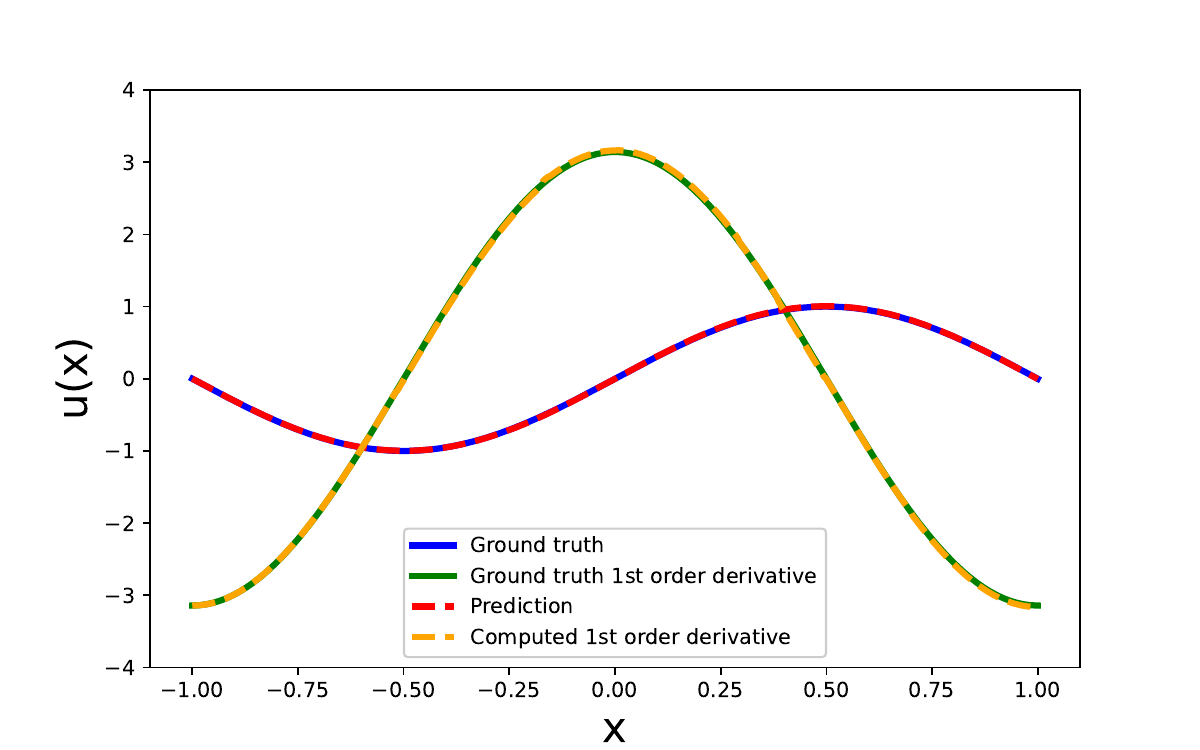}
		\caption{DSRN}
		\label{fig:compare-b}
	\end{subfigure}

	\begin{subfigure}[t]{0.45\textwidth}
		\centering
		\includegraphics[width=0.99\textwidth]{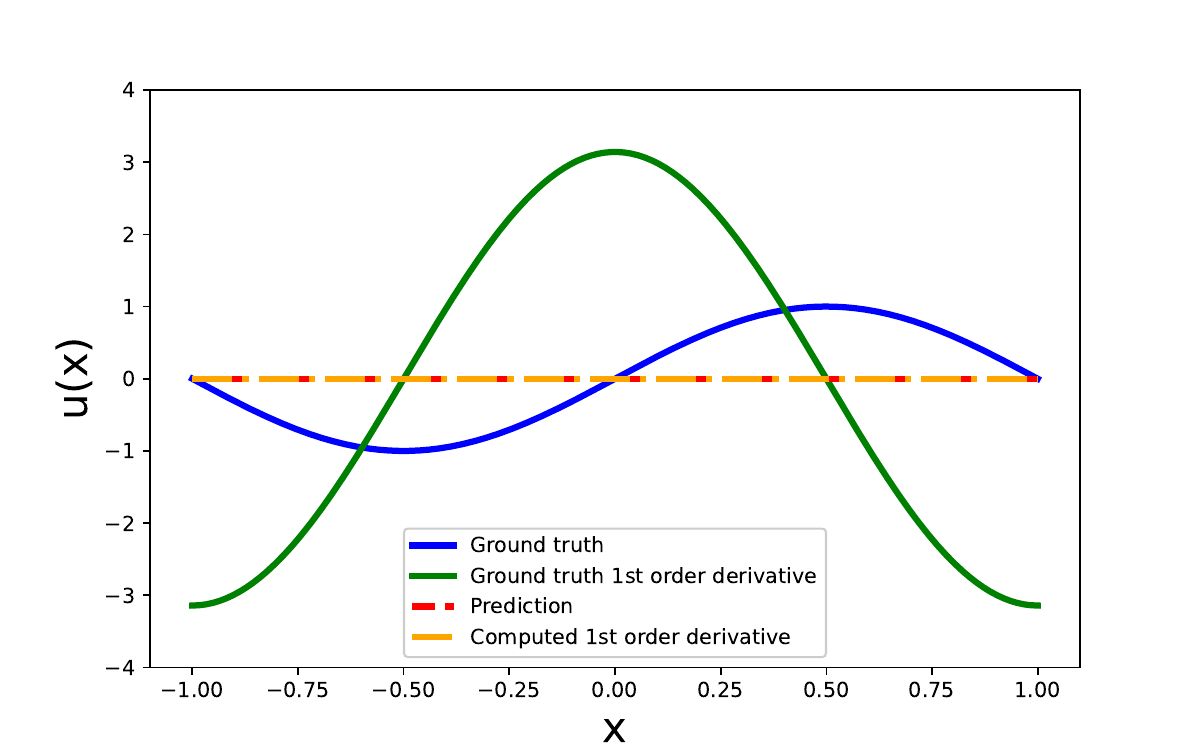}%V5
		\caption{Squared ReLU}
		\label{fig:compare-c}
	\end{subfigure}	
       \begin{subfigure}[t]{0.45\textwidth}
		\centering
		\includegraphics[width=0.99\textwidth]{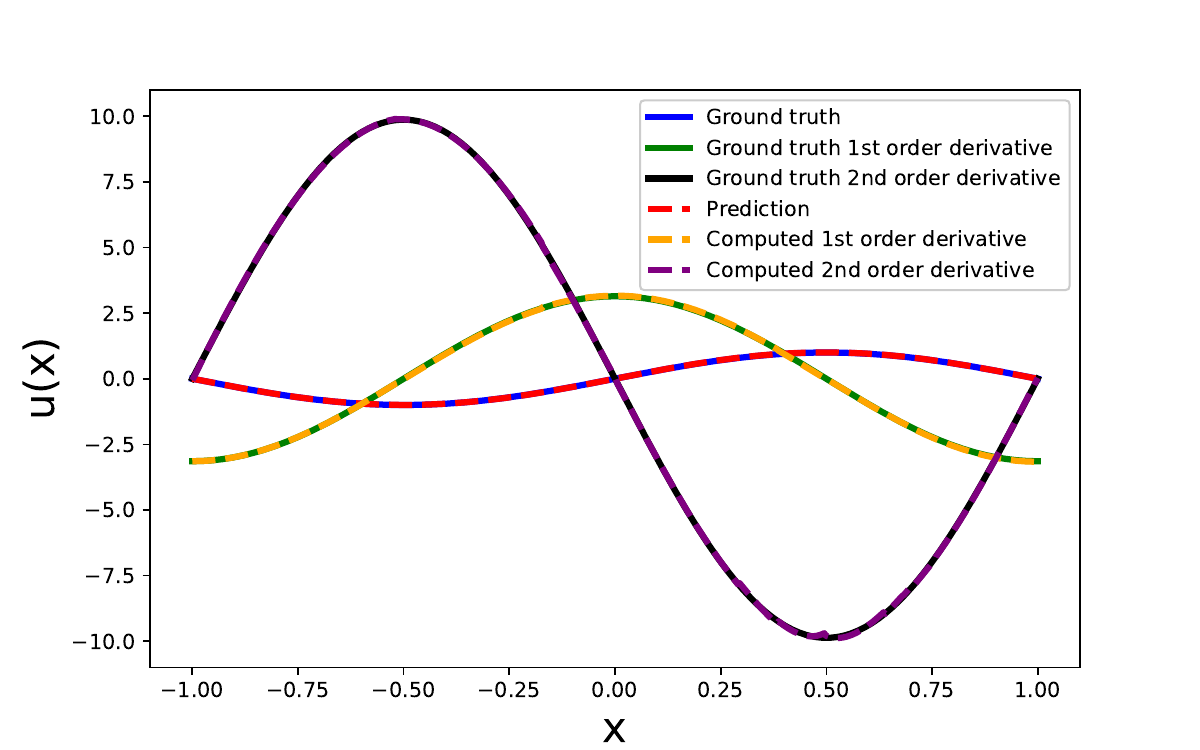}%V5
		\caption{DSRN}
		\label{fig:compare-d}
	\end{subfigure}	
	\caption{Comparison of prediction and the calculated first-order derivative between (a) ReLU, (b) DSRN and (c) Squared ReLU. (d) Visualization of the exact solution and the prediction, the calculated first-order, second-order derivative of DSRN. }	
	\label{fig:compare}
\end{figure}	

In Fig.~\ref{fig:compare},
We plot the prediction of the neural network with ReLU (Fig.~\ref{fig:compare-a}), DSRN (Fig.~\ref{fig:compare-b}) and the square of ReLU (Fig.~\ref{fig:compare-c}), respectively.
We also plot the first-order derivative $\frac{d}{d x} \mathcal{NN}_{\vtheta} $ using the autograd.
The exact solution is $u(x) = \sin(\pi x)$, which corresponds to the ground truth in the figure.
Fig.~\ref{fig:compare-a} shows that the networks only with ReLU activation function is not able to learn the solution with smooth first-order derivative.
However, in Fig.~\ref{fig:compare-b}, our method learns a solution with considerable smooth first-order derivative.
On the other hand, for the neural network with squared ReLU activation functions, we observed phenomena of gradient explosion and gradient vanishing during training, attributed to the depth of networks and the square operation.
Fig.~\ref{fig:compare-c} shows that the network using squared ReLU suffers from gradient vanishing and therefore fails to learn in the training.
These comparisons demonstrate that for deep fully connected networks, our DSRN method outperforms networks using only ReLU or squared ReLU.
We further visualize the second-order derivative of the DSRN outputs in Fig.~\ref{fig:compare-d}.
Fig.~\ref{fig:compare} shows that DSRN is able to approximate the smooth function with the square of ReLU activation function only in the last several layers in the neural network.
In Table \ref{tab:error}, we list the MSE between the neural network solutions and the exact solution: DSRN achieves more accurate approximation than neural networks with only ReLU or square of ReLU.
Note that the neural network with ReLU has zero-valued second-order derivative and thus a larger second-order derivative error.

\begin{table}[h!]
\centering
\caption{Mean squared error comparison.}
\begin{tabular}{|c|c|c|c|}
\hline
 & \textbf{ReLU} & \textbf{DSRN} & \textbf{Squared ReLU} \\
\hline
Solution error & $3 \times 10^{-4}$ & $\mathbf{1} \times \mathbf{10}^{-\mathbf{5}}$ & $0.5$ \\
\hline
Solution error + 2nd-derivative error & $48.61$ & $\mathbf{1.1} \times \mathbf{10}^{-\mathbf{3}}$ & $49.11$ \\
\hline
\end{tabular}\label{tab:error}
\end{table}

\subsection{Two-dimensional Diffusion equation}
\begin{figure}[h!]
\centering
        \begin{subfigure}[t]{0.48\textwidth}
		\centering
		\includegraphics[width=1.0\textwidth]{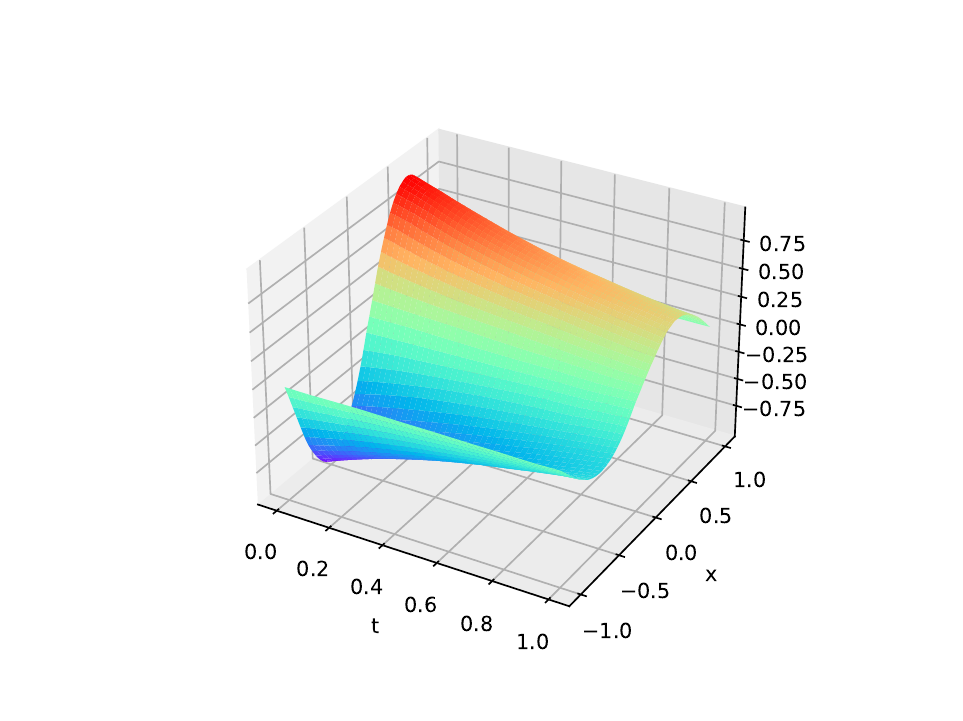}%V5
		\caption{DSRN output surface}
		\label{fig:compare-s-p}
	\end{subfigure}	
        \begin{subfigure}[t]{0.45\textwidth}
		\centering
		\includegraphics[width=0.90\textwidth]{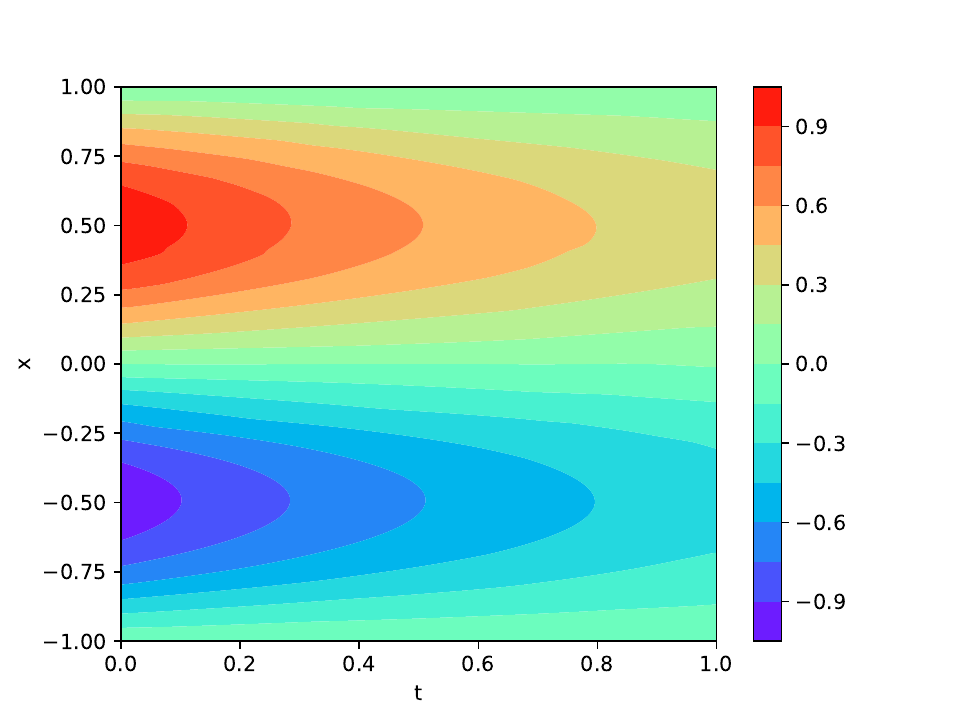}%V5
		\caption{DSRN output contour}
		\label{fig:compare-c-p}
	\end{subfigure}

	\begin{subfigure}[t]{0.48\textwidth}
		\centering
		\includegraphics[width=1.0\textwidth]{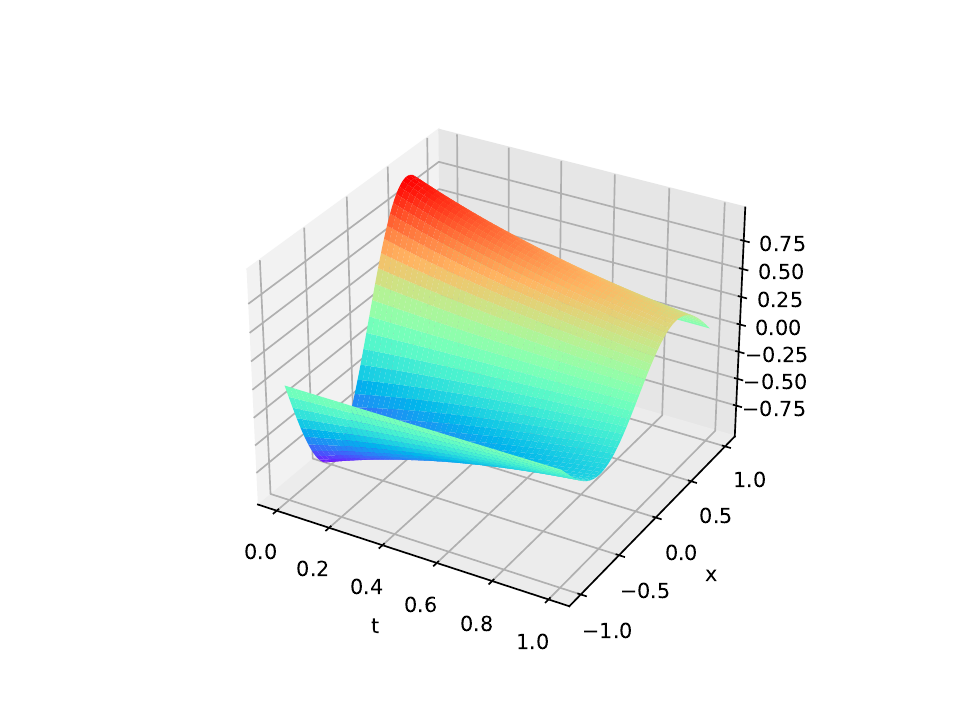}
		\caption{Exact solution surface}
		\label{fig:compare-s-t}
	\end{subfigure}
       \begin{subfigure}[t]{0.45\textwidth}
		\centering
		\includegraphics[width=0.90\textwidth]{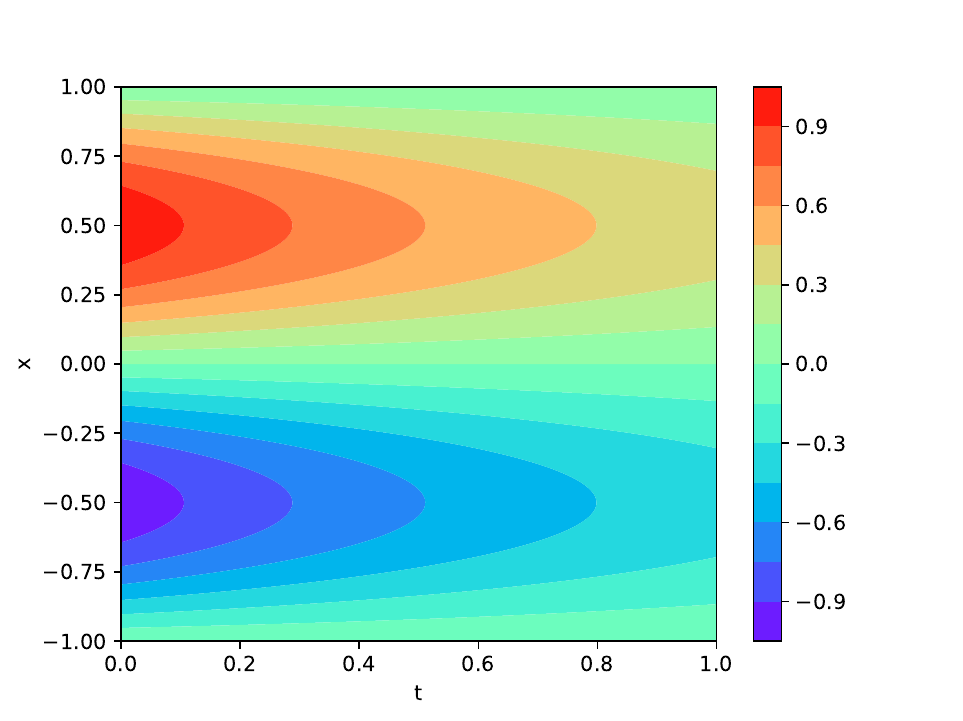}%V5
		\caption{Exact solution contour}
		\label{fig:compare-c-t}
	\end{subfigure}

	\caption{Comparison of surface and contour plots between DSRN predictions (a),(b) and the exact solution (c),(d).}	
	\label{fig:compare-2d}
\end{figure}

We solve the diffusion equation
\begin{equation}
\frac{\partial y}{\partial t} =\frac{\partial^2 y}{\partial x^2}-e^{-t}(\sin(\pi x)-\pi^2 \sin(\pi x)),
\quad 
x\in[-1,1], 
\quad
t\in[0,1],
\end{equation}
with initial and boundary conditions
$y(x,0)= \sin(\pi x), y(-1,t)=y(1,t)=0$.
The exact solution is $y(x,t)=e^{-t} \sin(\pi x)$.
We use a 12-layer fully-connected network $\mathcal{NN}_{\vtheta}(x,t)$ to learn the solution $y(x,t)$,
with the activation function of the last four hidden layers being the square of ReLU, and all other activation functions being ReLU.
Similarly, the loss function is \begin{footnotesize}
\[
\frac{1}{N} \sum_{i=1}^N 
\left[ 
\left(\frac{\partial \mathcal{NN}}{\partial t}(x_i,t_i) 
-\frac{\partial^2 \mathcal{NN}}{\partial x^2} (x_i,t_i)
+e^{-t_i}(sin(\pi x_i)-\pi^2 sin(\pi x_i)) \right)^2
+
\left(\mathcal{NN} (x_i, t_i) -  y_i \right)^2  
\right]\]    
\end{footnotesize}
where  $\left\{ (x_i,t_i, y_i),i=1,\dots,N | y_i=e^{-t_i} \sin(\pi x_i)\right\}$ are training data. 
In the training, we use the Adam optimizer with learning step $0.001$.
In Fig.~\ref{fig:compare-2d}, we plot surface and contours of the DSRN output and the exact solution, respectively.

\begin{figure}[h!]
\centering
\includegraphics[scale=0.50]{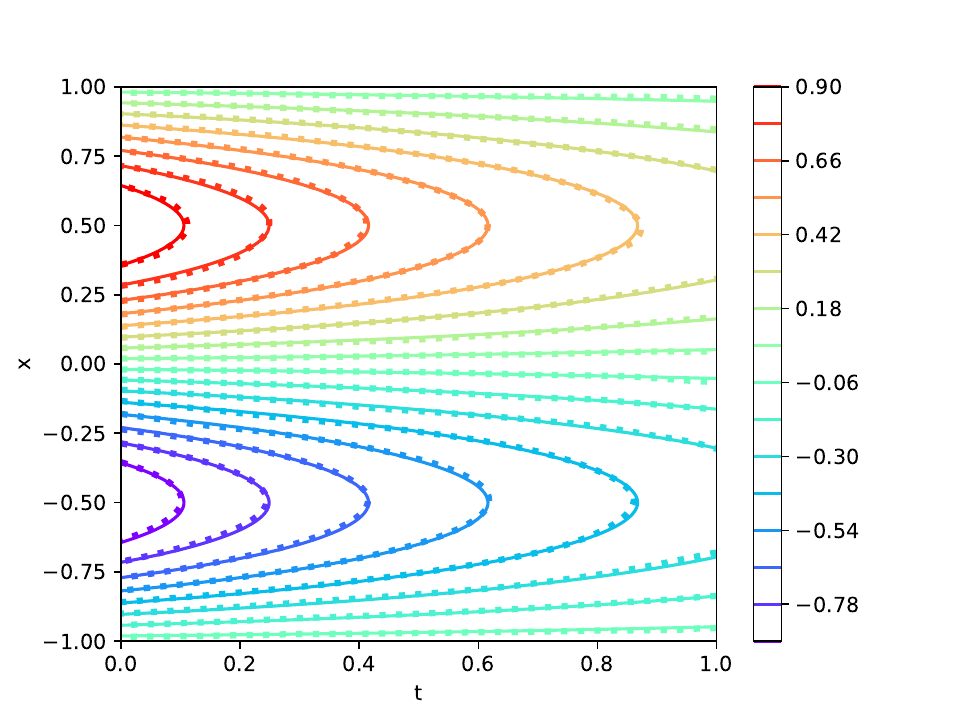}
\caption{Contour comparison: solid line for exact solution, dotted line for DSRN.}
\label{fig:compare-c-pt}
\end{figure}

For better comparison, we plot the contours of the DSRN output and the exact solution in one figure, see Fig.~\ref{fig:compare-c-pt}.
These comparisons demonstrate that DSRN accurately approximates the smooth solution.

\section{Conclusion}
           This paper introduces deep super ReLU networks (DSRNs) to approximate functions in Sobolev spaces measured by Sobolev norms, which ReLU DNNs cannot achieve. The DSRNs are constructed by adding several layers with the square of ReLU at the end of the neural networks to smooth them. While our DNNs are ReLU-based except for a few layers at the end, they retain the advantages of ReLU DNNs and are easy to train. Using other types of smooth activation functions instead of the square of ReLU may also achieve this goal, and we plan to explore these possibilities in future works. Additionally, the paper proves the optimality of DSRNs by estimating the VC-dimension of higher-order derivatives of DNNs and obtaining the generalization error in Sobolev space via the estimation of the pseudo-dimension of higher-order derivatives of DNNs. However, the method for estimating the VC-dimension of higher-order derivatives of DSRNs cannot be used for $\sigma_2$-NNs to obtain an optimal bound. Hence, a future research direction is to obtain optimal bounds for the VC-dimension of higher-order derivatives of $\sigma_2$-NNs.

\section*{Acknowledgments}
The work of H. Y. was partially supported by the US National
Science Foundation under award DMS-2244988, DMS-2206333, and the Office of Naval
Research Award N00014-23-1-2007. The work of Y. X. was supported by the Project of
Hetao Shenzhen-HKUST Innovation Cooperation Zone HZQB-KCZYB-2020083.

\bibliographystyle{elsarticle-num}
\bibliography{sample}
\end{document}